\documentclass[twoside]{article}
\usepackage{ecj,palatino,epsfig,latexsym,natbib,setspace}
\renewcommand{\cite}[2][]{\citep[{#1}]{{#2}}}

\usepackage{amsmath,amsfonts,amsthm}
\usepackage{url}
\usepackage{xcolor,graphicx}

\usepackage{ifthen}
\newcommand{\isdraft}{\boolean{true}} 
\renewcommand{\isdraft}{\boolean{false}} 
\ifthenelse{\isdraft}{
}{}
\newcommand{\markupdraft}[2]{
    \ifthenelse{\equal{#1}{display}}{#2}{}
    \ifthenelse{\equal{#1}{color}}{\color{#2}}{}
}
\newcommand{\mathnote}[1]{\markupdraft{display}{{\color{brown}\noindent[{\bf Math note}: #1]}}}
\newcommand{\nnotecolored}[3][]{\markupdraft{display}{{\color{#2}\noindent[{\bf Note #1}: #3]}}}
\newcommand{\newcolored}[3][]{{\markupdraft{color}{#2}#3}
    \ifthenelse{\equal{#1}{}}{}{\markupdraft{display}{{\color{yellow!70!black}[#1]}}}} 

\providecommand{\del}[2][]{{\markupdraft{display}{{\color{red!20!yellow}[rmed: "#2"[#1]]}}}} 
\providecommand{\new}[2][]{\newcolored[#1]{blue}{#2}}
\providecommand{\rem}[2][]{\nnotecolored[#1]{green}{#2}} 
   
\providecommand{\todo}[2][]{\markupdraft{display}{{\color{cyan}\noindent++TODO: #2 {\color{yellow}(#1)}++}}}
\ifthenelse{\isdraft}{}{\renewcommand{\markupdraft}[2]{}}
\newcommand{\niko}[1]{\rem[Niko]{\color{brown}#1}}
\newcommand{\anne}[1]{\rem[Anne]{\color{orange}#1}}
\newcommand{\alex}[1]{\rem[Alex]{\color{green}#1}}

\newcommand{\changeadd}[1]{{\color{red}#1}}
\newcommand{\changeaddp}[1]{{\color{blue}#1}} 
\newcommand{\changerem}[1]{{\color{green}#1}}
\newcommand{\changeremp}[1]{\changerem{#1}} 

\newcommand{\ad}[1]{\changeadd{#1}}
\newcommand{\re}[1]{\changerem{#1}}
\newcommand{\adp}[1]{\changeaddp{#1}}
\newcommand{\rep}[1]{\changerem{#1}}

\renewcommand{\changeaddp}[1]{{#1}}
\renewcommand{\changerem}[1]{{}}

\newcommand{\bs}[1]{\mathbf{#1}}

\newcommand{\R}{\mathbb{R}}
\newcommand{\N}{\ensuremath{\mathbb{N}}}
\newcommand{\E}{\mathbf{E}}
\newcommand{\borel}{\mathcal{B}}
\newcommand{\bplus}{\borel^{+}}
\newcommand{\1}{\mathbf{1}}
\newcommand{\mleb}{\mu_{Leb}}
\newcommand{\mr}{\rep{\mu_{Leb}}\adp{\mu_{\R_+^*}}}
\newcommand{\X}{\bs{X}}
\newcommand{\Xt}[1][t]{\X_{#1}}
\newcommand{\Xtt}{\X_{t+1}}
\newcommand{\w}{\bs{w}}
\newcommand{\W}{\mathcal{W}}
\newcommand{\Wt}[1][t]{\mathcal{W}_{#1}}

\newcommand{\x}{\bs{x}}
\newcommand{\y}{\bs{y}}
\newcommand{\Y}{\bs{Y}}
\newcommand{\Yti}[1][i]{\Y_t^{#1}}
\newcommand{\Ytilde}[1][i]{\bs{\tilde{Y}}_t^{#1}}

\newcommand{\Wti}[1][i]{\bs{W}_t^{#1}}

\newcommand{\G}{\mathcal{G}}
\newcommand{\Gtil}{\tilde{\G}}
\newcommand{\st}[1][t]{\sigma_{#1}}
\newcommand{\stt}{\sigma_{t+1}}
\newcommand{\stmod}[1][c]{\eta_{#1}^\star}
\newcommand{\Nti}[1][i]{\bs{N}_t^{#1}}
\newcommand{\Ntilde}[1][i]{\bs{\tilde{N}}_t^{#1}}
\newcommand{\Ntstar}[1][t]{\bs{\tilde{N}}_{#1}^\star}

\newcommand{\dt}[1][t]{\delta_{#1}}
\newcommand{\dmarkov}{\seq{\delta}}
\newcommand{\dpmarkov}{(\dt, \pt)_{t \in \N}}
\newcommand{\n}{\mathbf{n}}
\newcommand{\north}{\mathbf{n^{\perp}}}
\newcommand{\eone}{\mathbf{e}_{1}}
\newcommand{\etwo}{\mathbf{e}_{2}}
\newcommand{\ei}{\mathbf{e}_{i}}
\newcommand{\Nl}{\mathcal{N}}
\newcommand{\Nlun}{\mathcal{N}(0,1)}
\newcommand{\Nln}[1][n]{\mathcal{N}(\bs{0},\mathrm{Id}_{#1})}
\newcommand{\Nlambda}[1][\lambda]{\mathcal{N}_{#1:\lambda}}
\newcommand{\Ud}{\mathcal{U}}
\newcommand{\seq}[1]{\left( #1_t \right)_{t \in \N}}
\newcommand{\fdim}[1]{[#1]_1}

\newcommand{\chisq}[1][n-2]{\chi_{#1}^2}

\newcommand{\cdfnn}[1][\delta]{\tilde{F}_{#1}}
\newcommand{\cdfnninv}[1][\delta]{\tilde{F}_{#1}^{-1}}
\newcommand{\ptilde}[1][\delta]{p_{#1}}
\newcommand{\ptildeone}[1][\delta]{p_{1,#1}}

\newcommand{\pstar}[1][\delta]{p_{#1}^\star}
\newcommand{\pstarone}[1][\delta]{p_{1,#1}^\star}
\newcommand{\pstartwo}[1][\delta]{p_{2,#1}^\star}
\newcommand{\cdfnone}[1][\delta]{F_{1,\delta}}
\newcommand{\ddomain}{\R_+^{\changeaddp{*}}}
\newcommand{\doet}[1][t]{\dt[#1]}

\newcommand{\doedo}{\ddomain}
\newcommand{\doe}{\delta}
\newcommand{\equald}{\overset{d}{=}}
\newcommand{\pdfchi}{\varphi_\chi}
\newcommand{\p}{\bs{p}}
\newcommand{\pt}[1][t]{\p_{#1}}
\newcommand{\ptt}{\pt[t+1]}
\newcommand{\dttdomain}{\R_+^{\changeaddp{*}}}

\DeclareMathOperator{\argmax}{argmax}

\newtheorem{proposition}{Proposition}
\newtheorem{lemma}{Lemma}
\newtheorem{theorem}{Theorem}
\newtheorem{corollary}{Corollary}

\newenvironment{myproof}{\begin{proof}}{ \end{proof}} 

\begin{document}

\ecjHeader{x}{x}{xxx-xxx}{2015}{CSA-ES on a Linear \changeremp{Constraint}\changeaddp{Constrained} Problem}{A. Chotard, A. Auger, N. Hansen}
\title{\bf Markov Chain Analysis of Cumulative Step-size Adaptation on a Linear \changeremp{Constraint}\changeaddp{Constrained} Problem}  

\author{\name{\bf Alexandre Chotard} \hfill \addr{chotard@lri.fr}\\ 
        \addr{Univ. Paris-Sud, LRI, Rue Noetzlin, Bat 660, 91405 Orsay Cedex France}
\AND
       \name{\bf Anne Auger} \hfill \addr{auger@lri.fr}\\
       \addr{Inria, Univ. Paris-Sud, LRI, Rue Noetzlin, Bat 660, 91405 Orsay Cedex France}
\AND
       \name{\bf Nikolaus Hansen} \hfill \addr{hansen@lri.fr}\\
       \addr{Inria, Univ. Paris-Sud, LRI, Rue Noetzlin, Bat 660, 91405 Orsay Cedex France}
}

\maketitle

\begin{abstract}
This paper analyzes a $(1,\lambda)$-Evolution Strategy, a randomized comparison-based adaptive search algorithm, optimizing a linear function with a linear constraint. 
The algorithm uses resampling to handle the constraint. Two cases are investigated: first the case where the step-size is constant, and second the case where the step-size is adapted using cumulative step-size adaptation. 
We exhibit for each case a Markov chain describing the behaviour of the algorithm. 
Stability of the chain implies, by applying a law of large numbers, either convergence or divergence of the algorithm. Divergence is the desired behaviour. In the constant step-size case, we show stability of the Markov chain and prove the divergence of the algorithm. In the cumulative step-size adaptation case, we prove stability of the Markov chain in the simplified case where the cumulation parameter equals $1$, and discuss steps to obtain similar result\adp{s} for the full (default) algorithm where the cumulation parameter is smaller than $1$. The stability of the Markov chain allows us to deduce geometric divergence or convergence, depending on the dimension, constraint angle, population size and damping parameter, at a rate that we estimate.
Our results complement previous studies where stability was assumed.
\end{abstract}


\niko{I don't think that double braces are accepted by most editors. }

\niko{cleaning up the constraint/constrained change (and other trivial changes) would help readability with comments turned on. Generally I would suggest to clear up the changes that have proof read (in particular by the main author). }
\alex{
Cleaning up = making them not appear to us now? I kept them because I thought we should show to the reviewers what we deleted, so they can see better the changes we made.
}

\niko{I remain to find it utterly counterintuitive that green text is deleted text. It makes reading so much harder than it should be. Why didn't you just stick to the colors already defined? }
\alex{We were already using blue and yellow (not everything had been cleaned up then) and I wanted to have different colors that would mean changes for the reviewers from what they had seen. Then, I didn't not place special importance in using red and green, and I don't mind if you change it to something else.}

\section{Introduction}

Derivative Free Optimization (DFO) methods are  tailored for the optimization of numerical problems in a black-box context, where the objective function $f:\R^n \to \R$ is pictured as a black-box that \emph{solely} returns $f$ values (in particular no gradients are available).
\anne{Alexandre, I realized that Nesterov in a reference book calls also optimization with methods like BFGS, black-box optimization and this time the black box retunrs $f$ and $\nabla f$.}
\alex{I changed the text so that it talks solely about DFO.}\anne{It's not overly important but I guess it's good to mention the term "black-box optimization".}

Evolution Strategies (ES) are \emph{comparison-based} randomized DFO algorithms. At iteration $t$, solutions are sampled from a multivariate normal distribution centered in a vector $\Xt$. The candidate solutions are ranked according to $f$, and the updates of $\Xt$ and other parameters of the distribution (usually a step-size $\st$ and a covariance matrix) are performed using solely the ranking information given by the candidate solutions. Since ES do not directly use the function values of the new points, but only how the objective function $f$ ranks the different samples, they are invariant to the composition (to the left) of the objective function by a strictly increasing function $h: \R\to\R$.

This property and the black-box scenario make Evolution Strategies suited for a wide class of real-world problems, where constraints on the variables are often imposed. Different techniques for handling constraints in randomized algorithms have been proposed, see \rep{\cite{CoelloCoello:2008}}\adp{\cite{mezura2011constraint}} for a survey. For ES, common techniques are resampling, i.e.\ resample a solution until it lies in the feasible domain, repair of solutions that project unfeasible points onto the feasible domain \adp{\cite{arnold2011repair,arnold2013resamplingvsrepair}}, 
penalty methods where unfeasible solutions are penalised either by a quantity that depends on the distance to the constraint if this latter one can be computed (e.g.\ \adp{\cite{hansen2009tec,arnold2015towards}}\rep{\cite{hansen2009tec}} with adaptive penalty weights) or by the constraint value itself (e.g.\ stochastic ranking \cite{runarsson2000stochastic}) or methods inspired from multi-objective optimization (e.g.\ \cite{mezura2008constrained}).

In this paper we focus on the resampling method and study it on a simple \changeremp{constraint}\changeaddp{constrained} problem. More precisely, we study a $(1,\lambda)$-ES optimizing \emph{a linear function with a linear constraint} and resampling any infeasible solution until a feasible solution is sampled. The linear function models the situation where the current point is, relatively to the step-size, far from the optimum and ``solving'' this function means diverging. The linear constraint models being close to the constraint relatively to the step-size and far from other constraints. Due to the invariance of the algorithm to the composition of the objective function by a strictly increasing map, the linear function \adp{can}\rep{could} be composed by a function without derivative and with many discontinuities without any impact on our analysis.

The problem we address was studied previously  for different step-size adaptation mechanisms \adp{and different constraint handling methods}: with constant step-size, self-adaptation, and cumulative step-size adaptation\adp{, and the constraint being handled through resampling or repairing unfeasible solutions} \rep{\cite{arnold2011behaviour,arnold2012behaviour}}\adp{\cite{arnold2011behaviour,arnold2012behaviour,arnold2013resamplingvsrepair}}. The drawn conclusion is that when adapting the step-size the $(1,\lambda)$-ES fails to diverge unless some requirements on internal parameters of the algorithm are met. However, the approach followed in the aforementioned studies relies on finding simplified theoretical models to explain the behaviour of the algorithm: typically 
\adp{these models arise from} approximations (considering some random variables equal to their expected value, etc.) and assume mathematical properties like the existence of stationary distributions of underlying Markov chains \adp{without accompanied proof}. 

In contrast, our motivation is to study the \rep{real--i.e., not simplified--}algorithm \adp{without simplifications}\del{ these models arise from approximations}
and prove rigorously different mathematical properties of the algorithm allowing to deduce the exact behaviour of the algorithm, as well as to provide tools and methodology for such studies. Our theoretical studies need to be complemented by simulations of the convergence/divergence rates. The mathematical properties that we derive show that these numerical simulations converge fast. \adp{Our results are largely in agreement with the aforementioned studies of simplified models thereby backing up their validity.}

As for the step-size adaptation mechanism, our aim is to study the cumulative step-size adaptation (CSA) also called path-length control, default step-size  mechanism for the CMA-ES algorithm \cite{cmaes}. The mathematical object to study for this purpose is a discrete time, continuous state space Markov chain that is defined as the \changeremp{couple}\changeaddp{pair}:  evolution path and distance to the constraint \changeaddp{normalized by the step-size}. More precisely, stability properties like irreducibility and existence of a stationary distribution of this Markov chain need to be studied to deduce the geometric divergence of the CSA and have a rigorous mathematical framework to perform Monte Carlo simulations allowing to study the influence of different parameters of the algorithm. We start by illustrating in details the methodology on the simpler case where the step-size is constant. We show in this case that the distance to the constraint reaches a stationary distribution. This latter property was assumed in a previous study \cite{arnold2011behaviour}. We then prove that the algorithm  diverges at a constant speed. 
We then apply this approach to the case where the step-size is adapted using path length control. 
We show that in the special case where the cumulation parameter $c$ equals to $1$, 
the expected logarithmic step-size change, $\E\ln(\stt/\st)$, converges to a constant $r$, and the average logarithmic step-size change, $\ln(\st/\st[0])/t$, converges in probability to the same constant, which depends on parameters of the problem and of the algorithm. 
This implies geometric divergence (if $r>0$) or convergence (if $r<0$) at the rate $r$ \adp{for which estimations are provided}\rep{that we estimate}. 

This paper is organized as follows. In Section~\ref{sc:pr} we define the $(1,\lambda)$-ES using resampling and the problem. In Section~\ref{sec:preliminary} we provide some preliminary derivations on the distributions that come into play for the analysis. In Section~\ref{sc:cst} we analyze the constant step-size case. In Section~\ref{sc:csa} we analyze the cumulative step-size adaptation case. Finally we discuss our results and our methodology in Section~\ref{sc:discuss}.

A preliminary version of this paper appeared in the conference proceedings \cite{cah2014linearproblem}. The analysis of path-length control with cumulation parameter equal to $1$ is however fully new, as well as the discussion on how to analyze the case with cumulation parameter smaller than one. Also Figures~\ref{fg:csadvst}--\ref{fg:ccrit} are new as well as the convergence of the progress rate in Theorem~\ref{th:cstdiv}.

\subsection*{Notations}

Throughout this article, we denote by $\varphi$ the density function of the standard multivariate normal distribution \changeaddp{(the dimension being clarified within the context)}, and $\Phi$ the cumulative distribution function of a standard univariate normal distribution. The standard (unidimensional) normal distribution is denoted $\Nlun$, the ($n$-dimensional) multivariate  normal distribution with covariance matrix identity is denoted $\Nln$ and the $i^{\mathrm{th}}$ order statistic of $\lambda$ i.i.d. standard normal random variables is denoted $\Nlambda[i]$.  The uniform distribution on an interval $I$ is denoted $\Ud_I$. The set of natural numbers (including $0$) is denoted $\N$, and the set of real numbers $\R$. We denote $\R_+$ the set $\lbrace x \in \R | x \geq 0 \rbrace$, and for $A \subset \R^n$, the set $A^*$ denotes $A \backslash \lbrace\bs{0}\rbrace$  and $\1_{A}$ denotes the indicator function of $A$. For a topological space $\mathcal{X}$, $\mathcal{B}(\mathcal{X})$ denotes the Borel algebra of $\mathcal{X}$. We denote $\mleb$ the Lebesgue measure \adp{on $\R$, and for $A \subset \R$, $\mu_A$ denotes the trace measure $\mu_A: B \in \borel(\R) \mapsto \mleb(A \cap B)$}. For two vectors $\x \in \R^n$ and $\bs{y} \in \R^n$, we denote $[\x]_i$ the $i^{\mathrm{th}}$-coordinate of $\x$, and $\x .\y$ the scalar product of $\x$ and $\y$. Take $(a,b) \in \N^2$ with $a\leq b$, we denote $[a..b]$ the interval of integers between $a$ and $b$. \adp{The Gamma function is denoted by $\Gamma$.}  For $\X$ and $\Y$ two random vectors, we denote $\X \equald \Y$ if $\X$ and $\Y$ are equal in distribution. For $(X_t)_{t\in \N}$ a sequence of random variables and $X$ a random variable we denote $X_t \overset{a.s.}{\rightarrow} X$ if $X_t$ converges almost surely to $X$ and $X_t \overset{P}{\rightarrow} X$ if $X_t$ converges in probability to $X$. \adp{For $X$ a random variable and $\pi$ a probability measure, we denote $\E(X)$ the expected value of $X$, and $\E_\pi(X)$ the expected value of $X$ when $X$ has distribution $\pi$.}

\section{Problem statement and algorithm definition} \label{sc:pr}

\subsection{$(1,\lambda)$-ES with resampling}

In this paper, we study the behaviour of a $(1,\lambda)$-Evolution Strategy \emph{maximizing} a function $f$: $\R^n \rightarrow \R$, $\lambda \geq 2$, $n \geq 2$, with a constraint defined by a function $g : \R^n \rightarrow \R$ restricting the feasible space to $X_{\textrm{feasible}} = \lbrace \bs{x}\in \R^n | g(\bs{x}) \changerem{\geq}\changeadd{>} 0 \rbrace$. To handle the constraint, the algorithm resamples any unfeasible solution until a feasible solution is found.

From iteration $t \in \N$, given the vector $\Xt \in \R^{n}$ and step-size $\st \in \R_{+}^*$, the algorithm generates $\lambda$ new candidates:
\begin{equation}
\Yti = \Xt +\st \Nti \enspace,
\end{equation}
with $i \in [1 .. \lambda]$, and $(\Nti)_{i \in [1 .. \lambda]}$ i.i.d. standard multivariate normal random vectors. If a new sample $\Yti$ lies outside the feasible domain, that is $g(\Yti) \changeremp{<}\changeaddp{\leq} 0$, then it is resampled until it lies within the feasible domain. The first feasible $i^{\rm th}$ candidate solution is denoted $\Ytilde$ and the realization of the multivariate normal distribution giving $\Ytilde$ is $\Ntilde$, i.e.
\begin{equation}
\Ytilde = \Xt + \st \Ntilde
\end{equation}
 The vector $\Ntilde$ is called a feasible step. Note that $\Ntilde$ is not distributed as a multivariate normal distribution, further details on its distribution are given later on.

We define 
$\,\star = \argmax_{i \in [1 .. \lambda]} f( \Ytilde )$
as the index realizing the maximum objective function \adp{value}, and call $\Ntstar$ the selected step. The vector $\Xt$ is then updated as the solution realizing the maximum value of the objective function, i.e.
\begin{equation} \label{eq:Xtt}
\Xtt = \Ytilde[\star] = \Xt + \st \Ntstar \enspace.
\end{equation}

The step-size and other internal parameters are then adapted. We denote for the moment in a non specific manner the adaptation as
\begin{equation} \label{eq:sig}
\stt = \st \xi_t
\end{equation}
where $\xi_t$ is a random variable whose distribution is a function of the selected steps $(\Ntstar[i])_{i \leq t}$, $\Xt[0]$, $\st[0]$ and of internal parameters of the algorithm. \alex{I included in "internal parameters" $\pt[0]$ which is needed for $\xi_t$ but not introduced yet}We will define later on specific rules for this adaptation.

\subsection{Linear fitness function with linear constraint}

\begin{figure}
\begin{center}
\includegraphics[width=0.43\textwidth]{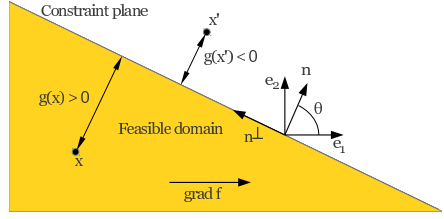}
\end{center}
\caption{Linear function with a linear constraint, in the plane generated by $\nabla f$ and $\n$, a normal vector to the constraint hyperplane with angle $\theta \in (0, \pi/2)$ with $\nabla f$. The point $\x$ is at distance $g(\x)$ from the constraint.}
\label{fg:problem}
\end{figure}

In this paper, we consider the case where $f$, the function that we optimize, and $g$, the constraint, are linear functions. W.l.o.g., we assume that $\| \nabla f\| = \| \nabla g \| = 1$. We denote $\n := - \nabla g $ a \changeremp{vector normal}\changeaddp{normal vector} to the constraint hyperplane. We choose an orthonormal Euclidean coordinate system with basis $(\ei)_{i \in [1 ..n]}$ with its origin located on the constraint hyperplane where $\eone$ is equal to the gradient $\nabla f$, hence 
\begin{equation}\label{eq:f}
f(\x) = [\x]_{1}
\end{equation}
and the vector $\etwo$ lives in the plane generated by $\nabla f$ and $\n$ and is such that the angle between $\etwo$ and $\n$ is positive. We define $\theta$ the angle between $\nabla f$ and $\n$, and restrict our study to $\theta \in (0, \pi/2)$. The function $g$ can be seen as a signed distance to the linear constraint as
\begin{equation} \label{eq:g}
g(\x) = \x.\nabla g = -\x.\n  =  -[\x]_1 \cos \theta -[\x]_2 \sin \theta \enspace .
\end{equation}
A point is feasible if and only if $g(\x) \changeremp{\geq}\changeaddp{>} 0$ (see Figure~\ref{fg:problem}). Overall the problem reads
\begin{equation}\label{eq:pbdef}
\begin{split}
{\rm maximize}\,\,\,\, f(\x) = [\x]_{1}\,\,\,\, {\rm subject~to~}  \\ g(\x) = -[\x]_1 \cos \theta -[\x]_2 \sin \theta \changeremp{\geq}\changeaddp{>} 0 \enspace.
\end{split} 
\end{equation}

Although $\Ntilde$ and $\Ntstar$ are in $\R^n$, due to the choice of the coordinate system and  the independence of the sequence $([\Nti]_k)_{k \in [1..n]}$, only the two first coordinates of these vectors are affected by the resampling implied by $g$ and the selection according to $f$. Therefore $[\Ntstar]_k \sim \Nlun$ for $k \in [3..n]$. With an abuse of notations, the vector $\Ntilde$ will  denote the 2-dimensional vector $([\Ntilde]_1,[\Ntilde]_2)$, likewise $\Ntstar$ will also denote the 2-dimensional vector $([\Ntstar]_{1},[\Ntstar]_{2})$, and $\n$ will denote the 2-dimensional vector $(\cos \theta, \sin \theta)$. The coordinate system will also be used as $(\eone, \etwo)$  only.

Following \cite{arnold2011behaviour,arnold2012behaviour,arnold2008behaviour}, we denote the normalized signed distance to the constraint as $\dt$, that is 
\begin{equation}\label{eq:dt}
\dt = \frac{g(\Xt)}{\st} \enspace .
\end{equation}

We initialize the algorithm by choosing $\Xt[0] = -\n$ and $\st[0] = 1$, which implies that $\dt[0] = 1$.

\section{Preliminary results and definitions}\label{sec:preliminary}

Throughout this section we derive the probability density functions of the random vectors $\Ntilde$ and $\Ntstar$ and give a definition of $\Ntilde$ and of $\Ntstar$ as a function of $\dt$ and of an i.i.d. sequence of random vectors.

\subsection{Feasible steps} \label{sc:fea}

The random vector $\Ntilde$, the $i^{\rm th}$ feasible step, is distributed as the standard multivariate normal distribution truncated by the constraint, as stated in the following lemma.
\begin{lemma} \label{lm:fea}
Let a $(1,\lambda)$-ES with resampling optimize a function $f$ under a constraint function $g$. If $g$ is a linear form determined by a vector $\n$ as in \eqref{eq:g}, then the distribution of the feasible step $\Ntilde$ only depends on the normalized distance to the constraint $\dt$ and its density given that $\dt$ equals $\delta$ reads
\begin{equation} \label{eq:p}
\ptilde \left( \x  \right) = \frac{ \varphi(\x) \1_{\dttdomain}\left( \delta - \x.\n \right)}{\Phi(\delta)} \enspace .
\end{equation}
\end{lemma}

\begin{myproof}
A solution $\Yti$ is feasible if and only if $g(\Yti) \changeremp{\geq}\changeaddp{>} 0$, which is equivalent to $-(\Xt+\st \Nti).\n \changeremp{\geq}\changeaddp{>} 0 $. Hence dividing by $\st$, a solution is feasible if and only if $\dt  = - \Xt.\n/\st \changeremp{\geq}\changeaddp{>} \Nti.\n $.
Since a standard multivariate normal distribution is rotational invariant, $\Nti.\n$ follows a standard (unidimensional) normal distribution. Hence the probability that a solution $\Yti$ or a step $\Nti$  is feasible is given by 
\begin{equation*} 
\Pr( \Nlun \changeremp{\leq}\changeaddp{<} \dt ) = \Phi\left(\dt\right) \enspace.
\end{equation*}
Therefore the probability density function of the random variable $\Ntilde.\n$ for $\dt = \delta$ is $x \mapsto \varphi(x) \1_{\dttdomain}(\delta - x)/\Phi(\delta)$. For any vector $\north$ orthogonal to $\n$ the random variable $\Ntilde.\north$ was not affected by the resampling and is therefore still distributed as a standard (unidimensional) normal distribution. With a change of variables \alex{From the previous point, we have that 
$$
\Pr(\Ntilde.\n \leq \tilde{x} \textrm{ and }\Ntilde.\north \leq  \tilde{y} | \delta_t = \delta) = \int_{\R} \int_{\R} \frac{\varphi(\tilde{x}) \varphi(\tilde{y}) \1_{\R_+^*}(\delta-\tilde{x})}{\Phi(\delta)} d\tilde{x} d\tilde{y} \enspace .
$$
With the change of variables $\tilde{x} = x\cos \theta  + y\sin \theta $ and $\tilde{y} = -x\sin \theta  + y\cos \theta $, and since $\Ntilde.\n \leq \tilde{x}$ and $\Ntilde.\north \leq \tilde{y}$ is equivalent to $[\Ntilde]_1 \leq x$ and $[\Ntilde]_2 \leq y$, we obtain the joint distribution of Eq.~\eqref{eq:p}.} using the fact that the standard multivariate normal distribution is rotational invariant we obtain the joint distribution of Eq.~\eqref{eq:p}. 
\end{myproof}

Then the marginal density function $\ptildeone$ of $[\Ntilde]_{1}$ can be computed by integrating Eq.~\eqref{eq:p} over $[\x]_{2}$ and reads
\begin{equation} \label{eq:p1}
\ptildeone\left(x \right) = \varphi \left(x\right) \frac{\Phi\left(\frac{\delta - x \cos \theta}{\sin \theta}\right)}{\Phi\left(\delta \right)} \enspace ,
\end{equation}
(see \cite[Eq.~4]{arnold2011behaviour} for details) and we denote $\cdfnone$ its cumulative distribution function.

It will be important in the sequel to be able to express the vector $\Ntilde$ as a function of $\dt$ and of a \emph{finite} number of random samples. Hence we give an alternative way to sample $\Ntilde$ rather than the resampling technique that involves an unbounded number of samples.

\begin{lemma} \label{lm:Ntilde}
Let a $(1,\lambda)$-ES with resampling optimize a function $f$ under a constraint function $g$, where $g$ is a linear form determined by a vector $\n$ as in \eqref{eq:g}.
Let the feasible step $\Ntilde$ be the random vector described in Lemma~\ref{lm:fea} and $\bs{Q}$ be the 2-dimensional rotation matrix of angle $\theta$. Then 
\begin{equation}\label{def:Ntilde}
\Ntilde ~\equald~ \cdfnninv[\dt](U_{t}^{i}) \n + \mathcal{N}_{t}^{i} \north = \bs{Q}^{-1} \begin{pmatrix} \cdfnn[\dt]^{-1}(U_{t}^{i}) \\ \mathcal{N}_{t}^{i} \end{pmatrix}
\end{equation}
where $\cdfnninv[\dt]$ denotes the generalized inverse of the cumulative distribution of $\Ntilde . \n$ \footnote{The generalized inverse of $\cdfnn$ is $\cdfnninv[\delta](y) := \inf_{x\in\R}\lbrace \cdfnn[\delta](x) \geq y \rbrace$.}, $U_{t}^{i} \sim \Ud_{[0,1]}$, $\mathcal{N}_{t}^{i} \sim \Nlun$ with $(U_{t}^{i})_{i \in [1..\lambda], t \in \N}$ i.i.d. and $(\mathcal{N}_{t}^{i})_{i \in [1..\lambda], t \in \N}$ i.i.d. random variables. 
\end{lemma}

\begin{myproof}%
We define a new coordinate system $(\n, \north)$ (see Figure~\ref{fg:problem}). It is the image of $(\eone, \etwo)$ by  $\bs{Q}$. In the new basis $(\n, \north)$, only the coordinate along $\n$ is affected by the resampling. Hence the random variable $\Ntilde.\n$  follows a truncated normal distribution with cumulative distribution function $\cdfnn[\dt]$ equal to $\min(1,\Phi(x)/\Phi(\dt))$, while the random variable $\Ntilde.\north$ follows an independent standard normal distribution, hence $\Ntilde \equald (\Ntilde.\n)\n + \mathcal{N}_{t}^{i} \n^{\perp}$. Using the fact that if a random variable has a cumulative distribution $F$, then for $F^{-1}$ the generalized inverse of $F$, $F^{-1}(U)$ with $U \sim \Ud_{[0,1]}$ has the same distribution as this random variable, we get that $\cdfnn[\dt]^{-1}(U_{t}^{i}) \equald \Ntilde.\n$, so we obtain Eq.~\eqref{def:Ntilde}.
\end{myproof}

We now extend our study to the selected step $\Ntstar$.

\subsection{Selected step}

The selected step $\Ntstar$ is chosen among the different feasible steps $(\Ntilde)_{i \in [1..\lambda]}$ to maximize the function $f$, and has the density described in the following lemma.
\begin{lemma} \label{lm:pstar}
Let a $(1,\lambda)$-ES with resampling optimize the problem \eqref{eq:pbdef}. Then the distribution of the selected step $\Ntstar$ only depends on the normalized distance to the constraint $\dt$ and its density given that $\dt$ equals $\delta$ reads 
\begin{align} \label{eq:pstar}
\pstar \! \left(\x \right) \! &= \! \lambda \ptilde\left(\x \right) \cdfnone([\x]_1 )^{\lambda -1} \enspace , \\
 &= \! \lambda \frac{\varphi(\x) \1_{\dttdomain}(\delta - \x.\n)}{\Phi(\delta)} \!\left( \! \int_{-\infty}^{[\x]_1} \!\!\!\!\!\! \varphi(u) \frac{\Phi(\frac{\delta - u\cos\theta}{\sin \theta})}{\Phi(\delta)} \mathrm{d}u \!\right)^{\lambda-1} \nonumber
\end{align}
where $\ptilde$ is the density of $\Ntilde$ given that $\dt = \delta$ given in Eq.~\eqref{eq:p} and $\cdfnone$ the cumulative distribution function of $[\Ntilde]_1$ whose density is given in Eq.~\eqref{eq:p1} and $\n$ the vector $(\cos\theta, \sin\theta)$.
\end{lemma}

\begin{myproof}
The function $f$ being linear, the rankings on $(\Ntilde)_{i\in [1..\lambda]}$ correspond to the order statistic on $([\Ntilde]_1)_{i\in [1..\lambda]}$. If we look at the joint cumulative distribution $F_\delta^\star$ of $\Ntstar$
\begin{align*}
F_\delta^\star(x,y) &= \Pr \left( [\Ntstar]_1 \leq x, [\Ntstar]_2 \leq y\right) \\
& = \sum_{i=1}^\lambda \Pr \left( \!  \Ntilde[i] \leq \left( \!\!\begin{array}{c} x \\ y \end{array} \!\! \right), [\Ntilde[j]]_1 < [\Ntilde[i]]_1 \text{ for } j \neq i  \right)
\end{align*}
by summing disjoints events. 
The vectors $(\Ntilde)_{i \in [1..\lambda]}$ being independent and identically distributed 
\begin{align*}
F_\delta^\star(x,y) &= \lambda \Pr \left( \Ntilde[1] \leq \left( \! \begin{array}{c} x \\ y \end{array} \! \right), [\Ntilde[j]]_1 < [\Ntilde[1]]_1 \text{ for } j \neq 1 \right) \\
& = \lambda \int_{-\infty}^x\int_{-\infty}^y p_\delta (u,v) \prod_{j=2}^\lambda \Pr([\Ntilde[j]]_1 < u) \mathrm{d}v \mathrm{d}u \\
& = \lambda \int_{-\infty}^x\int_{-\infty}^y p_\delta (u,v) \cdfnone(u)^{\lambda-1} \mathrm{d}v \mathrm{d}u \enspace .
\end{align*}
Deriving $F_\delta^\star$ on $x$ and $y$ yields the density of $\Ntstar$ of Eq.~\eqref{eq:pstar}.
\end{myproof}

We may now obtain the marginal of $[\Ntstar]_1$ and $[\Ntstar]_2$. 
\begin{corollary}
 Let a $(1,\lambda)$-ES with resampling optimize the problem \eqref{eq:pbdef}. Then the marginal distribution of $[\Ntstar]_1$ only depends on $\dt$ and its density given that $\dt$ equals $\delta$ reads
\begin{align} \label{eq:pstar1}
\pstarone \left(x \right) &= \lambda \ptildeone ( x ) \cdfnone(x )^{\lambda -1} \enspace , \\
&= \lambda \varphi(x) \frac{\Phi\left(\frac{\delta - x\cos\theta}{\sin\theta}\right)}{\Phi(\delta)} \cdfnone(x )^{\lambda -1} \enspace , \nonumber 
\end{align}
and the same holds for $[\Ntstar]_2$ whose marginal density reads
\begin{align} \label{eq:pstar2}
\pstartwo \left(y \right) &= \lambda \frac{\varphi(y)}{\Phi(\delta)} \int_{-\infty}^{\frac{\delta - y \sin \theta}{\cos \theta}} \varphi(u) \cdfnone(u)^{\lambda-1} \mathrm{d}u \enspace.
\end{align}
\end{corollary}

\begin{myproof}
Integrating Eq.~\eqref{eq:pstar} directly yields Eq.~\eqref{eq:pstar1}.

The conditional density function of $[\Ntstar]_2$ is 
$$
\pstartwo (y | [\Ntstar]_1 = x) = \frac{\pstar( (x,y) )}{ \pstarone(x )} \enspace .
$$
As $\pstartwo(y ) = \int_{\R} \pstartwo (y | [\Ntstar]_1 = x) \pstarone (x) \mathrm{d}x$, using the previous equation with Eq.~\eqref{eq:pstar} gives that $\pstartwo(y ) = \int_{\R} \lambda \ptilde((x,y) ) \cdfnone (x)^{\lambda-1} \mathrm{d}x$, which with Eq.~\eqref{eq:p} gives 
\begin{equation*}
 \pstartwo (y ) = \lambda \frac{\varphi(y)}{\Phi(\delta)}\int_{\R} \! \varphi(x)\1_{\R_+^{\changeadd{*}}} \! \left(\delta - \left( \!\! \begin{array}{c} x \\ y \end{array} \!\! \right) .\n \right) \cdfnone(x)^{\lambda-1} \mathrm{d}x.
\end{equation*}
The condition $\delta - x\cos\theta - y\sin\theta \rep{\geq}\adp{>} 0$ is equivalent to $x \rep{\leq}\adp{<} (\delta - y\sin\theta)/\cos\theta$, hence Eq.~\eqref{eq:pstar2} holds.
\end{myproof}

We will need in the next sections an expression of the random vector $\Ntstar$ as a function of $\dt$ and a random vector composed of a \emph{finite} number of i.i.d. random variables. To do so, using notations of Lemma~\ref{lm:Ntilde}, we define the function $\Gtil : \ddomain\times([0,1]\times\R) \rightarrow \R^2$ as
\begin{equation} \label{eq:Gtilde}
\Gtil(\delta, \bs{w}) = \bs{Q}^{-1} \left( \begin{array}{c} \cdfnn^{-1}\left( [\bs{w}]_1 \right) \\ {[\bs{w}]_2} \end{array} \right) \enspace .
\end{equation}
According to Lemma~\ref{lm:Ntilde}, given that $U \sim \Ud_{[0,1]}$ and $\Nl \sim \Nlun$, $( \cdfnn^{-1}(U) , \Nl )$ (resp. $\Gtil(\delta,(U,\Nl))$) is distributed as the resampled step $\Ntilde$ in the coordinate system $(\n, \n^\perp)$ (resp. $(\bs{e}_1, \bs{e}_2)$).
 Finally, let $(\bs{w}_i)_{i \in [1..\lambda]} \in ([0,1]\times \R)^\lambda$ and let $\G : \ddomain \times([0,1]\times\R)^{\lambda} \rightarrow \R^2$ be the function defined as
\begin{equation} \label{eq:G}
\G(\delta, (\bs{w}_i)_{i \in [1..\lambda]}) = \underset{\bs{N} \in \left\lbrace \Gtil(\delta, \bs{w}_i) | i \in [1..\lambda] \right\rbrace}{\argmax} f(\bs{N}) \enspace .
\end{equation}
As shown in the following proposition, given that $\Wti \sim (\Ud_{[0,1]}, \Nlun)$ and $\Wt = (\Wti)_{i \in [1 .. \lambda]}$, the function $\G(\dt, \Wt)$ is distributed as the selected step $\Ntstar$.

\begin{proposition} \label{pr:g}
Let a $(1,\lambda)$-ES with resampling optimize the problem defined in Eq.~\eqref{eq:pbdef}, and let $(\Wti)_{i \in [1 .. \lambda], t \in \N}$ be an i.i.d. sequence of random vectors with $\Wti \sim (\Ud_{[0,1]}, \Nlun)$, and $\Wt = (\Wti)_{i \in [1 .. \lambda]}$. Then
\begin{equation}\label{eq:NtstarG}
\Ntstar  ~\equald~  \G(\dt, \Wt) \enspace,
\end{equation}
where the function $\G$ is defined in Eq.~\eqref{eq:G}.
\end{proposition}

\begin{myproof} 
Since $f$ is a linear function $f(\Ytilde) = f(\Xt) + \st f(\Ntilde)$, so $f(\Ytilde) \leq f(\Ytilde[j])$ is equivalent to $f(\Ntilde) \leq f(\Ntilde[j])$. Hence $\star = \argmax_{i \in [1..\lambda]} f(\Ntilde)$ and therefore $\Ntstar = \argmax_{\bs{N} \in \lbrace \Ntilde | i \in [1..\lambda] \rbrace } f(\bs{N})$. From Lemma~\ref{lm:Ntilde} and Eq.~\eqref{eq:Gtilde}, $\Ntilde \equald  \Gtil(\dt, \Wti)$, so $\Ntstar \equald \argmax_{\bs{N} \in \lbrace \Gtil(\dt, \Wti) | i \in [1..\lambda] \rbrace } f(\bs{N})$, which from \eqref{eq:G} is $\G(\dt, \Wt)$.
\end{myproof}

\section{Constant step-size case} \label{sc:cst}

We illustrate in this section our methodology \changeremp{analysis} on the simple case where the step-size is constantly equal to $\sigma$ and prove that  $(\Xt)_{t\in\N}$ diverges in probability at constant speed and that the progress rate $\varphi^* := \E([\Xtt]_1 - [\Xt]_1) = \sigma\E([\Ntstar]_1)$ (see \citealt[Eq.~2]{arnold2011behaviour}) converges to a strictly positive constant  (Theorem~\ref{th:cstdiv}). The analysis of the CSA is then a generalization of the results presented here, with  more technical results to derive. Note that the progress rate definition coincides with the fitness gain, i.e. $\varphi^*=\E(f(\Xtt) - f(\Xt))$.

As suggested in \cite{arnold2011behaviour}, the sequence $\dmarkov$ plays a central role for the analysis, and we will show that it admits a stationary measure. We first prove that this sequence is a homogeneous Markov chain.
\begin{proposition}
Consider the $(1,\lambda)$-ES with resampling and with constant step-size $\sigma$ optimizing the \changeremp{constraint}\changeaddp{constrained} problem~\eqref{eq:pbdef}. 
Then the sequence $\dt=g(\Xt)/\sigma$ is a homogeneous Markov chain on $\ddomain$ and
\begin{equation} \label{eq:cstmarkov}
\dt[t+1] = \dt - \Ntstar . \n ~\equald~ \dt - \G(\dt, \Wt).\n \enspace, 
\end{equation}
where $\G$ is the function defined in~\eqref{eq:G} and $(\Wt)_{t\in\N} = (\Wti)_{i \in [1..\lambda],t\in\N}$ is an i.i.d. sequence with $\Wti \sim (\Ud_{[0,1]}, \Nlun)$ for all $(i,t) \in [1..\lambda]\times \N$.
\end{proposition}
\begin{myproof}
It follows from the definition of $\dt$ that $\dt[t+1] = \frac{g\left(\Xtt\right)}{\st[t+1]} = \frac{-\left(\Xt + \sigma\Ntstar\right).\n}{\sigma} = \dt - \Ntstar . \n$, and in Proposition~\ref{pr:g} we state that $\Ntstar ~\equald~ \G(\dt, \Wt)$. Since $\dt[t+1]$ has the same distribution as a time independent function of $\dt$ and of $\Wt$ where $(\Wt)_{t \in \N}$ are i.i.d., it is a homogeneous Markov chain.
\end{myproof}

The Markov Chain $\dmarkov$ comes into play for investigating the divergence of $f(\Xt)=[\Xt]_{1}$. Indeed, we can express $\frac{[\Xt - \Xt[0]]_1}{t}$ in the following manner:
\begin{align}
\frac{[\Xt - \Xt[0]]_1}{t} & = \frac{1}{t} \sum_{k=0}^{t-1} \left( [\X_{k+1}]_{1} - [\X_{k}]_{1} \right) \nonumber
 \\ \label{eq:xrec} &= 
\frac{\sigma}{t}\sum_{k=0}^{t-1} [\Ntstar[k]]_{1} ~\equald~ \frac{\sigma}{t}\sum_{k=0}^{t-1} [\G(\dt[k],\W_{k})]_1 \enspace .
\end{align}
The latter term suggests the use of a Law of Large Numbers (LLN) to prove the convergence of $\frac{[\Xt - \Xt[0]]_1}{t}$ which will in turn imply--\changeadd{-}if the limit is positive\changeadd{-}--the divergence of $[\Xt]_1$ at a constant rate. Sufficient conditions on a Markov chain to be able to apply the LLN include the existence of an invariant probability measure $\pi$. The limit term is then expressed as an expectation over the stationary distribution. More precisely, assume the LLN can be applied, the following limit will hold
\begin{align}
\frac{[\Xt - \Xt[0]]_1}{t} \overset{a.s.}{\underset{t \to \infty}{\longrightarrow}} \sigma \int_{\ddomain} \E  \left([\G(\delta,\W) ]_1\right) \pi(d \delta) \label{eq:xlln} \enspace.
\end{align}

If the Markov chain $\dmarkov$ is also $V$-ergodic with $|\E([G(\delta,\W)]_1)| \leq V(\delta)$ then the progress rate converges to the same limit.
\begin{equation}
\E([\Xtt]_1 -[\Xt]_1) \underset{t \to +\infty}{\longrightarrow} \sigma \int_{\ddomain} \E  \left([\G(\delta,\W) ]_1\right) \pi(d \delta) \label{eq:progressesp} \enspace.
\end{equation}
We prove formally these two equations in Theorem~\ref{th:cstdiv}.

The invariant measure $\pi$ is also underlying  the study carried out in \cite[Section~4]{arnold2011behaviour} where more precisely it is stated: {\it ``Assuming for now that the mutation strength $\sigma$ is held constant, when the algorithm is iterated, the distribution of $\delta$-values tends to a stationary limit distribution.''}. We will now provide a formal proof that indeed $\dmarkov$ admits a stationary limit distribution $\pi$, as well as prove some other useful properties that will allow us in the end to conclude to the divergence of $([\Xt]_1)_{t \in \N}$.

\subsection{Study of the stability of $\dmarkov$}

We study in this section the stability of $\dmarkov$. We first derive its transition kernel $P(\delta,A) := \Pr( \dt[t+1] \in A | \dt[t] = \delta) $ for all $\delta \in \ddomain$ and $A \in \borel(\ddomain)$. Since
$\Pr( \dt[t+1] \in A | \dt = \delta)  = \Pr(\delta_{t} - \Ntstar.\n \in A | \dt = \delta ) \enspace,
$
\begin{equation}\label{trans:kernel}
P(\delta,A) = \int_{\R^2} \1_A \left( \delta - \bs{u}.\n \right) \pstar \left( \bs{u}\right) \mathrm{d}\bs{u}
\end{equation}
where $\pstar$ is the density of $\Ntstar$ given in \eqref{eq:pstar}. For $t \in \N^*$, the $t$-steps transition  kernel $P^t$ is defined by $P^t(\delta, A) := \Pr(\dt \in A | \dt[0] = \delta)$.

From the transition kernel, we will now derive the first properties on the Markov chain $\dmarkov$. First of all we investigate the so-called $\psi$-irreducible property.

A Markov chain $(\doet)_{t\in\N}$ on a state space $\doedo$ is  \emph{$\psi$-irreducible} if there exists a non-trivial measure $\psi$ such that for all sets $A \in \borel(\doedo)$ with $\psi(A)>0$ and for all $\doe \in \doedo$, there exists $t \in \N^*$ such that $P^t(\doe,A)>0$. We denote $\bplus(\doedo)$ the set of Borel sets of $\doedo$ with strictly positive $\psi$-measure.

We also need the notion of \emph{small sets} \changeaddp{and \emph{petite sets}}. A set $C \in \borel(\doedo)$ is called a small set if there exists $m \in \N^*$ and a non trivial measure $\nu_m$ such that for all sets $A \in \borel(\doedo)$ and all $\delta \in C$
\begin{equation}
P^m(\doe,A) \geq \nu_m(A) \enspace .
\end{equation}
\changeaddp{A set $C \in \borel(\doedo)$ is called a petite set if there exists a probability measure $\alpha$ on $\N$ and a non trivial measure $\mu_\alpha$ such that for all sets $A \in \borel(\doedo)$ and all $\delta \in C$
\begin{equation}
K_\alpha(x, A) := \sum_{m \in \N} P^m(\x, A) \alpha(m) \geq \mu_\alpha(A) \enspace .
\end{equation}
A small set is therefore also a petite set. As we will see further, the existence of a small set combined with a control of the Markov chain chain outside of the small set allows to deduce powerful stability properties of the Markov chain.}
If there exists a $\nu_1$-small set $C$ such that $\nu_1(C) > 0$ then the Markov chain is said \emph{strongly aperiodic}.

\begin{proposition} \label{pr:cstirreducible}
Consider a $(1,\lambda)$-ES with resampling and with constant step-size optimizing the \changeremp{constraint}\changeaddp{constrained} problem \eqref{eq:pbdef}  and let $\dmarkov$ be the Markov chain exhibited in \eqref{eq:cstmarkov}. Then $\dmarkov$ is $\mr$-irreducible, strongly aperiodic, and compact sets \changeaddp{of $\R_+^*$ and sets of the form $(0,M]$ with $M > 0$} are small sets.
\anne{so $C$ is a compact of $\R^*_+$ ? Could you specify. Note that,  given that the state space is $\R^*_+$, I do not see that it makes sense to talk about compact of $\R$ (because the small set needs to be an element of $\borel(\doedo)$). But hence we also have that $(0,M]$ is a compact of $ \R^*_+$ or Am I Wrong here (I am confused) and thus you should rather state the proposition and "compact sets of $ \R^*_+$ and in particular sets of the form $(0,M]$ with $M > 0$ are small sets.}
\end{proposition}

\begin{myproof}
\adp{Take $\delta \in \R_+^*$ and $A \in \borel(\R_+^*)$.}
Using Eq.~\eqref{trans:kernel} and Eq.~\eqref{eq:pstar} the transition kernel can be written
\begin{equation*}
P(\delta, A) \! = \! \lambda \! \int_{\R^2} \!\! \1_{A} (\delta - \left( \!\! \begin{array}{c} x \\ y \end{array} \!\! \right).\n ) \frac{\varphi(x)\varphi(y)}{\Phi(\delta)} \cdfnone( x )^{\lambda-1} \mathrm{d}y\mathrm{d}x \enspace .
\end{equation*}
We remove $\delta$ from the indicator function by a substitution of variables  $u = \delta - x\cos\theta -y\sin\theta$, and $v = x\sin\theta -y\cos\theta$. As this substitution is the composition of a rotation and a translation the determinant of its Jacobian matrix is $1$. We denote $h_\delta : (u,v) \mapsto (\delta-u)\cos\theta + v\sin\theta$, $h_\delta^\perp : (u,v) \mapsto (\delta-u)\sin\theta - v\cos\theta$ and $g (\delta,u,v) \mapsto \lambda \varphi(h_\delta(u,v))\varphi(h_\delta^\perp(u,v))/\Phi(\delta)\cdfnone(h_\delta(u,v))^{ \lambda -1}$. Then $x = h_\delta(u,v)$, $y = h_\delta^\perp(u,v)$ and
\begin{equation}
\label{eq:ker}
P(\delta, A) = \int_{\R}  \int_\R 1_{A}(u) g(\delta,u,v)  \mathrm{d}v \mathrm{d}u \! \enspace .
\end{equation}
For all $\delta,u,v$ the function $g(\delta,u,v)$ is strictly positive hence for all $A$ with $\mr(A) > 0$, $P(\delta,A) > 0$. Hence $\dmarkov$ is irreducible with respect to the Lebesgue measure.

In addition, the function $(\delta,u,v) \mapsto g(\delta,u,v)$ is continuous as the composition of continuous functions (the continuity of $\delta \mapsto F_{1,\delta}(x)$ for all $x$ coming from the  dominated convergence theorem). Given a compact $C$ \adp{of $ \R^*_+$,} we hence know that there exists $g_C > 0$ such that for all $(\delta,u,v) \in C\times[0,1]^2$, $g(\delta,u,v) \geq g_C > 0$.
Hence for all $\delta \in C$,
$$
P(\delta,A) \geq \underbrace{g_C  \mr(A \cap [0,1]) }_{:=\nu_C(A)} \enspace.
$$
The measure $\nu_C$ being non-trivial, the previous equation shows that compact sets \adp{of $ \R^*_+$,} are small and that for $C$ a compact such that $\mr(C \cap [0,1]) > 0$, we have $\nu_C(C) > 0$ hence the chain is strongly aperiodic. \changeaddp{Note also that since $\lim_{\delta \to 0} g(\delta, u, v) > 0$, the same reasoning holds for $(0,M]$ instead of $C$ (where $M > 0$). Hence the set $(0,M]$ is also a small set.}\anne{mhum, so something needs to be adapted if $C$ is right from the start a compact set of $ \R^*_+$.}
\end{myproof}

The application of the LLN for a $\psi$-irreducible Markov chain $(\doet)_{t \in \N}$ on a state space $\doedo$ requires the existence of an \emph{invariant measure} $\pi$, that is satisfying for all $A \in \borel(\doedo)$
\begin{equation} \label{eq:positive}
\pi(A) = \int_{\doedo} P(\doe, A) \pi(\mathrm{d}\doe)  \enspace .
\end{equation}
If a Markov chain admits an invariant probability measure then the Markov chain is called positive.

A typical assumption to apply the LLN is positivity and Harris-recurrence.
A $\psi$-irreducible chain $(\doet)_{t \in \N}$ on a state space $\doedo$ is \emph{Harris-recurrent} if for all sets $A \in \bplus(\doedo)$ and for all $\doe \in \doedo$, $\Pr ( \eta_A = \infty | \doe_0 = \doe ) = 1$ where $\eta_A$ is the occupation time of A, i.e. $\eta_A = \sum_{t=1}^\infty 1_A(\doet)$.
We will show that the Markov chain $\dmarkov$ is positive and Harris-recurrent by using so-called Foster-Lyapunov drift conditions: define the \emph{drift} operator for a positive function $V$ as 
\begin{equation} \label{eq:driftop}
\Delta V(\delta) = \E[ V(\dt[t+1]) | \dt=\delta ] - V(\delta)\enspace.
\end{equation}
Drift conditions translate that outside a small set, the drift operator is negative. We will show a drift condition for V-geometric ergodicity where given a function $f \geq 1$, a positive and Harris-recurrent chain $(\doet)_{t \in \N}$ with invariant measure $\pi$ is called \emph{$f$-geometrically ergodic} if $\pi(f) \changeaddp{:= \int_{\R} f(\doe)\pi(\textrm{d}\doe)} < \infty$ and there exists $r_f> 1$ such that
\begin{equation} \label{eq:geo}
\sum_{t \in \N} r_f^t \| P^t(\doe, \cdot) - \pi \|_f < \infty \enspace , \forall \doe  \in  \doedo \enspace ,
\end{equation}
where for $\nu$ a signed measure $\| \nu \|_f$ denotes $\sup_{g : |g| \leq f} | \int_{\doedo} g(x)\nu(\textrm{d}x) |$.

To prove the $V$-geometric ergodicity, we will prove that there exists a small set $C$, constants $b \in \R$, $\epsilon \in \R_+^*$ and a function $V \geq 1$ finite for at least some $\doe_0 \in \doedo$ such that for all $\doe \in \doedo$
\begin{equation} \label{eq:V4}
\Delta V(\delta) \leq -\epsilon V(\delta) + b \1_{C}(\delta) \enspace .
\end{equation}
If the Markov chain $\dmarkov$ is $\psi$-irreducible and aperiodic, this drift condition implies that the chain is $V$-geometrically ergodic \cite[Theorem 15.0.1]{markovtheory}\footnote{The condition $\pi(V) < \infty$ is given by \cite[Theorem 14.0.1]{markovtheory}.}  as well as positive and Harris-recurrent\footnote{The function $V$ of \eqref{eq:V4} is unbounded off small sets \cite[Lemma 15.2.2]{markovtheory} with \cite[Proposition 5.5.7]{markovtheory}, hence with \cite[Theorem 9.1.8]{markovtheory} the Markov chain is Harris-recurrent.}.

Because \rep{compacts}\adp{sets of the form $(0,M]$ with $M> 0$} are small sets and drift conditions investigate the negativity outside a small set, we need to study the chain for $\delta$ large. The following lemma is a technical lemma studying the limit of $\E(\exp(\G(\delta, \W).\n))$  for $\delta$ to infinity.

\begin{lemma} \label{lm:farNtstar}
Consider the $(1,\lambda)$-ES with resampling optimizing the \changeremp{constraint}\changeaddp{constrained} problem \eqref{eq:pbdef}, and let $\G$ be the function defined in \eqref{eq:G}.
We denote $K$ and $\bar{K}$ the random variables $\exp(\G(\delta, \W).(a,b))$ and $\exp(a|[\G(\delta,\W)]_1| + b|[\G(\delta,\W)]_2|)$. For $\W \sim (\Ud_{[0,1]}, \Nl(0,1))^\lambda$ and any $(a,b)\in \R^2$, $\lim_{\delta \rightarrow +\infty} \E(K) = \E(\exp(a\Nlambda))\E(\exp(b \Nlun)) < \infty$ and $\lim_{\delta \rightarrow +\infty} \E(\bar{K})  < \infty$
\end{lemma}
For the proof see the appendix. We are now ready to prove a drift condition for geometric ergodicity.

\begin{proposition} \label{pr:csterg}
Consider a $(1,\lambda)$-ES with resampling and with constant step-size optimizing the \changeremp{constraint}\changeaddp{constrained} problem \eqref{eq:pbdef} and let $\dmarkov$ be the Markov chain exhibited in \eqref{eq:cstmarkov}. The Markov chain $\dmarkov$ is $V$-geometrically ergodic with $V : \delta \mapsto \exp(\alpha \delta)$ for $\alpha > 0$ small enough, and is Harris-recurrent and positive with invariant probability measure $\pi$.
\end{proposition}

\begin{myproof}
Take the function $V : \delta \mapsto \exp(\alpha\delta)$, then 
\begin{align*}
\Delta V (\delta) &= \E\left(\exp\left(\alpha(\delta - \G(\delta,\W).\n )\right)\right) - \exp(\alpha \delta ) \\
\frac{\Delta V}{V}(\delta)	&= \E\left(\exp\left( - \alpha\G(\delta,\W).\n \right)\right) - 1 \enspace . 
\end{align*}
With Lemma~\ref{lm:farNtstar} we obtain that
$$
\underset{\delta\rightarrow + \infty}{\lim} \E\left(\exp\left( - \alpha\G(\delta,\W).\n \right)\right) =  \E\left( \exp (-\alpha \Nlambda \cos\theta) \right) \E(\exp(-\alpha \Nlun \sin \theta)) < \infty \enspace .
$$
As the right hand side of the previous equation is finite we can invert integral with series with Fubini's theorem, so with Taylor series
\begin{equation*}
\underset{\delta\rightarrow + \infty}{\lim} \!\! \E\left(\exp\left( - \alpha\G(\delta,\W).\n \right)\right) = \left( \! \sum_{i \in \N} \! \frac{\left(-\alpha \cos\theta\right)^i \E \! \left(\Nlambda ^i\right)}{i!} \!\! \right) \!\! \left( \! \sum_{i \in \N} \! \frac{\left(-\alpha \sin\theta \right)^i \E \! \left(\Nlun^i\right)}{i!} \!\! \right) \!\! \enspace ,
\end{equation*}
which in turns yields
\begin{align*}
\lim_{\delta \rightarrow + \infty} \!\! \frac{\Delta V}{V}(\delta) &= \! \left(1 - \alpha\E(\Nlambda)\cos\theta + o(\alpha)\right)\left(1 + \! o(\alpha)\right) \! - \! 1 \\
&= -\alpha\E(\Nlambda)\cos\theta + o(\alpha) \enspace .
\end{align*}
Since for $\lambda \geq 2$, $\E(\Nlambda) > 0$, for $\alpha >0$ and small enough we get $\lim_{\delta \rightarrow + \infty} \frac{\Delta V}{V}(\delta) < -\epsilon < 0$. Hence there exists $\epsilon > 0$, $M>0$ and $b\in \R$ such that 
\begin{equation*}
\Delta V(\delta) \leq -\epsilon V(\delta) + b \1_{\changeremp{[}\changeaddp{(}0,M]}(\delta) \enspace .
\end{equation*}

According to Proposition~\ref{pr:cstirreducible}, $\changeremp{[}\changeaddp{(}0,M]$ is a small set, hence it is petite \cite[Proposition 5.5.3]{markovtheory}. Furthermore $\dmarkov$ is a $\psi$-irreducible aperiodic Markov chain so $\dmarkov$ satisfies the conditions of Theorem~15.0.1 from \cite{markovtheory}, which with Lemma~15.2.2, Theorem~9.1.8 and Theorem~14.0.1 of \cite{markovtheory} proves the proposition.
\end{myproof}

We now proved rigorously the existence (and unicity) of an invariant measure $\pi$ for the Markov chain $\dmarkov$, which provides the so-called steady state behaviour in \cite[Section 4]{arnold2011behaviour}. 
As the Markov chain $\dmarkov$ is positive and Harris-recurrent we may now apply a Law of Large Numbers \cite[Theorem 17.1.7]{markovtheory} in Eq~\eqref{eq:xrec} to obtain the divergence of $f(\Xt)$ and an exact expression of the divergence rate.

\begin{theorem} \label{th:cstdiv} 
Consider a $(1,\lambda)$-ES with resampling and with constant step-size optimizing the \changeremp{constraint}\changeaddp{constrained} problem \eqref{eq:pbdef} and let $\dmarkov$ be the Markov chain exhibited in \eqref{eq:cstmarkov}. 
The sequence  $(\fdim{\Xt})_{t \in \N}$ diverges in probability to $+\infty$ at constant speed, that is
\begin{align} \label{eq:cstdiv}
\frac{\fdim{\Xt - \Xt[0]}}{t} & \overset{P}{\underset{t \rightarrow +\infty}{\longrightarrow}} \sigma  \E_{\pi\otimes\mu_\W} \left( [\G\left(\delta, \W \right)]_1\right)  > 0 
\end{align}
and the expected progress satisfies
\begin{align} \label{eq:expdiv}
 \varphi^* = \E\left( \fdim{\Xtt - \Xt} \right) & \underset{t \rightarrow +\infty}{\longrightarrow} \sigma  \E_{\pi\otimes\mu_\W} \left( [\G\left(\delta, \W \right)]_1\right)  > 0 \enspace ,
\end{align}
where $\varphi^*$ is the progress rate defined in \cite[Eq.~(2)]{arnold2011behaviour}, $\G$ is defined in \eqref{eq:G}, $\W = (\bs{W}^i)_{i \in [1..\lambda]}$ with $(\bs{W}^i)_{i \in [1..\lambda]}$ an i.i.d. sequence such that $\bs{W}^i \sim (\Ud_{[0,1]}, \Nlun)$, $\pi$ is the stationary measure of $\dmarkov$ whose existence is proven in Proposition~\ref{pr:csterg} and $\mu_\W$ is the probability measure of $\W$.
\end{theorem}

\begin{myproof} From Proposition~\ref{pr:csterg} the Markov chain $\dmarkov$ is  Harris-recurrent and positive, and since $(\Wt)_{t \in \N}$ is i.i.d., the chain $(\dt, \Wt)$ is also Harris-recurrent and positive with invariant probability measure $\pi \times \mu_\W$, so to apply the Law of Large Numbers \cite[Theorem 17.0.1]{markovtheory} to $[\G]_1$ we only need $[\G]_1$ to be $\pi\otimes\mu_\W$-integrable.

With Fubini-Tonelli's theorem $\E_{\pi\otimes \mu_\W}(|[\G(\delta,\W )]_1|) $ equals to $ \E_{\pi} (\E_{\mu_\W}(|[\G(\delta,\W )]_1| ))$. As $\delta\geq 0$, we have $\Phi(\delta) \geq \Phi(0) = 1/2$, and for all $x\in\R$  as  $\Phi(x) \leq 1$, $\cdfnone(x) \leq 1$ and $\varphi(x) \leq \exp(-x^2/2)$ with Eq.~\eqref{eq:pstar1} we obtain that $|x| \pstarone(x) \leq 2 \lambda |x| \exp(-x^2/2)$ so the function $x \mapsto |x| \pstarone(x)$ is integrable. Hence for all $\delta \in \R_+$, $\E_{\mu_\W}(|[\G(\delta,\W )]_1| )$ is finite. Using the dominated convergence theorem, the function $\delta \mapsto \cdfnone(x)$ is continuous, hence so is $\delta \mapsto \pstarone(x)$. From \eqref{eq:pstar1} $|x|\pstarone(x) \leq 2 \lambda |x| \varphi(x)$,  which is integrable, so the dominated convergence theorem implies that the function $\delta \mapsto \E_{\mu_\W}(|[\G(\delta,\W]_1|)$ is continuous. Finally, using Lemma~\ref{lm:farNtstar} with Jensen's inequality shows that $\lim_{\delta \rightarrow +\infty} \E_{\mu_\W}(|[\G(\delta, \W)]_1|)$ is finite. Therefore the function  $\delta \mapsto \E_{\mu_\W}(|[\G(\delta,\W]_1|)$ is bounded by a constant $M \in \R_+$. As $\pi$ is a probability measure $\E_{\pi} (\E_{\mu_\W}(|[\G(\delta,\W )]_1| )) \leq M < \infty$, meaning $[\G]_1$ is $\pi\otimes\mu_\W$-integrable. Hence we may apply the LLN on Eq.~\eqref{eq:xrec}
\begin{equation*}
\frac{\sigma}{t} \sum_{k=0}^{t-1} [\G(\dt[k],\Wt[k])]_1 \overset{a.s.}{\underset{t \rightarrow + \infty}{\longrightarrow}} \sigma \E_{\pi\otimes\mu_\W}\left([\G(\delta,\W )]_1\right) < \infty \enspace .
\end{equation*}
The equality in distribution in \eqref{eq:xrec} allows us to deduce the convergence in probability of the left hand side of \eqref{eq:xrec} to the right hand side of the previous equation.

From \eqref{eq:xrec} $[\Xtt - \Xt]_1 \equald \sigma\G(\dt, \Wt)$ so $\E([\Xtt - \Xt]_1 | \Xt[0] = \x) =  \sigma \E(\G(\dt, \Wt) | \dt[0] = \x/\sigma)$. As $\G$ is integrable with Fubini's theorem $\E(\G(\dt, \Wt) | \dt[0] = \x/\sigma) = \int_{\ddomain}  \E_{\mu_\W}(\G( \y, \W)) P^t(\x/\sigma, d\y) $, so $ \E(\G(\dt, \Wt) | \dt[0] = \x/\sigma) - \E_{\pi \times \mu_\W}(\G(\delta, \W)) = \int_{\ddomain} \E_{\mu_\W}(\G( \y, \W)) (P^t(\x/\sigma, d\y) - \pi(d\y))  $. According to Proposition~\ref{pr:csterg} $\dmarkov$ is $V$-geometrically ergodic with $V : \delta \mapsto \exp(\alpha \delta)$, so there exists $M_\delta$ and $r > 1$ such that $\|P^t(\delta, \cdot) - \pi \|_V \leq M_\delta r^{-t} $. We showed that the function $\delta \mapsto \E(|[\G(\delta, \W)]_1|)$ is bounded, so since $V(\delta) \geq 1$ for all $\delta \in \ddomain$ and $\lim_{\delta \rightarrow +\infty} V(\delta) = + \infty$, there exists $k$ such that $\E_{\mu_\W}(|[\G(\delta, \W)]_1|) \leq k V(\delta)$ for all $\delta$. Hence $| \int \E_{\mu_\W}(|[\G( x, \W)]_1|) (P^t(\delta, dx) - \pi(dx)) | \leq  k\|P^t(\delta, \cdot) - \pi \|_V \leq k M_\delta r^{-t} $. And therefore $|\E(\G(\dt, \Wt) | \dt[0] = \x/\sigma) - \E_{\pi \times \mu_\W}(\G(\delta, \W))| \leq k M_\delta r^{-t}$ which converges to $0$ when $t$ goes to infinity.

As the measure $\pi$ is an invariant measure for the Markov chain $\dmarkov$, using \eqref{eq:cstmarkov}, $\E_{\pi\otimes\mu_\W}(\delta) = \E_{\pi\otimes\mu_\W}(\delta-\G(\delta,\W).\n)$, hence $\E_{\pi\otimes\mu_\W}(\G(\delta,\W).\n) = 0$ and thus
\begin{equation*}
\E_{\pi\otimes\mu_\W}\left([ \G(\delta, \W)]_1 \right) = - \tan \theta \E_{\pi\otimes\mu_\W}\left([\G(\delta, \W)]_2\right) \enspace .
\end{equation*}
We see from Eq.~\eqref{eq:pstar2} that for $y > 0$, $\pstartwo(y ) < \pstartwo(-y)$  hence the expected value $\E_{\pi\otimes\mu_\W}([\G(\delta, \W)]_2)$ is strictly negative. With the previous equation it implies that $\E_{\pi\otimes\mu_\W}([\G(\delta, \W)]_1)$ is strictly positive.

\end{myproof}

We showed rigorously the divergence of $[\Xt]_1$ and gave an exact expression of the divergence rate, and that the progress rate $\varphi^*$ converges to the same rate. The fact that the chain $\dmarkov$ is $V$-geometrically ergodic gives that there exists a constant $r>1$ such that $\sum_t r^t \|P^t(\delta,\cdot) - \pi\|_V < \infty$. This implies that the distribution $\pi$ can be simulated efficiently by a Monte Carlo simulation allowing to have precise estimations of the divergence rate of $[\Xt]_1$.

\begin{figure}\centering
\includegraphics[width=0.65\textwidth, trim=0 0ex 0 5ex 0, clip]{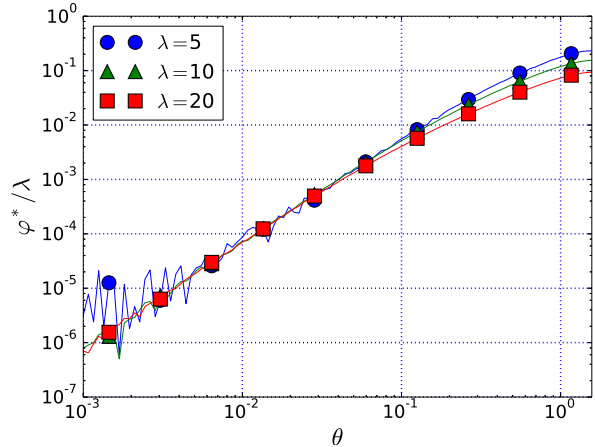}
\caption{Normalized progress rate $\varphi^{\rep{\star}\adp{*}} = \E(f(\Xtt) - f(\Xt))$ divided by  $\lambda $ for the $(1,\lambda)$-ES with constant step-size $\sigma = 1$ and resampling, plotted against the constraint angle $\theta$, for $\lambda \in \lbrace 5, 10, 20 \rbrace$.}
\label{fg:cstpvst}
\end{figure}

A Monte Carlo simulation of the divergence rate in the right hand side of \eqref{eq:cstdiv} and \eqref{eq:expdiv} and for $10^6$ time steps gives the progress rate of \cite{arnold2011behaviour} $\varphi^{\re{\star}\ad{*}} =  \E([\Xtt - \Xt]_1)$,  which once normalized by $\sigma$ and $\lambda$ yields Fig.~\ref{fg:cstpvst}. We normalize per $\lambda$ as in evolution strategies the cost of the algorithm is assumed to be the number of $f$-calls. We see that for small values of $\theta$, the normalized serial progress rate assumes roughly $\varphi^{\re{\star}\ad{*}}/\lambda \approx \theta^2$. Only for larger constraint angles the serial progress rate depends on $\lambda$ where smaller $\lambda$ are preferable. 

\begin{figure}\centering
\includegraphics[width=0.65\textwidth, trim=0 0ex 0 5ex 0, clip]{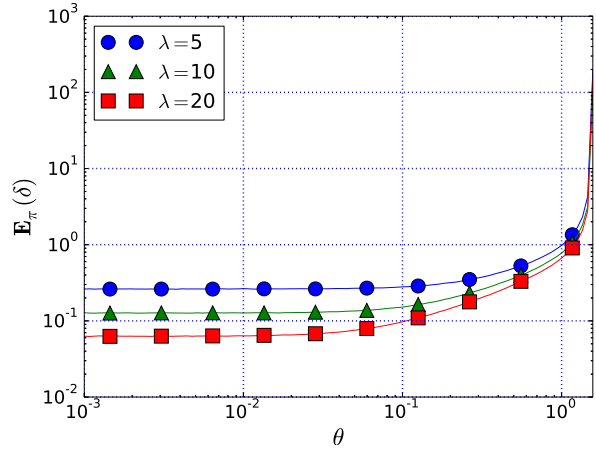}
\caption{Average normalized distance $\delta$ from the constraint for the $(1,\lambda)$-ES with constant step-size and resampling plotted against the constraint angle $\theta$ for $\lambda \in \lbrace 5, 10, 20 \rbrace$.}
\label{fg:cstdvst}
\end{figure}

Fig.~\ref{fg:cstdvst} is obtained through simulations of the Markov chain $\dmarkov$ defined in Eq.~\eqref{eq:cstmarkov} for $10^6$ time steps where the values of $\dmarkov$ are averaged over time. We see that when $\theta \rightarrow \pi/2$ then $\E_\pi ( \delta) \rightarrow + \infty$ since the selection does not attract $\Xt$ towards the constraint anymore. With a larger population size the algorithm is closer to the constraint, as better samples are more likely to be found close to the constraint.

\section{Cumulative Step size Adaptation} \label{sc:csa}

In this section we apply the techniques introduced in the previous section to the case where the step-size is adapted using Cumulative Step-size Adaptation. This technique was studied on sphere functions \cite{arnold2004performance} and on ridge functions \cite{arnold2008step}.

In CSA, the step-size is adapted using a path $\pt$, vector of $\mathbb{R}^n$, that sums up the different selected steps $\Ntstar$ with a discount factor. More precisely the evolution path $\pt \in \R^n$ is defined by $\pt[0] \sim \Nln$ and
\begin{equation} \label{eq:path}
\ptt = (1-c)\pt + \sqrt{c(2-c)}\Ntstar \enspace .
\end{equation} The variable $c \in (0,1]$ is called the cumulation parameter, and determines the "memory" of the evolution path, with the importance of a step $\Ntstar[0]$ decreasing in $(1-c)^t$. The backward time horizon is consequently about $1/c$. The coefficients in Eq~\eqref{eq:path} are chosen such that if $\pt$ follows a standard normal distribution, and if $f$ ranks uniformly randomly the different samples $(\Ntilde)_{i \in [1\changerem{,}\changeadd{..}\lambda]}$ \changeaddp{and that these samples are normally distributed}, then $\pt[t+1]$ will also follow a standard normal distribution independently of the value of $c$.

The length of the evolution path is compared to the expected length of a Gaussian vector (that corresponds to the expected length under random selection) (see \cite{cmaes}). To simplify the analysis we study here a modified version of CSA introduced in \cite{arnold2002noisy} where the squared length of the evolution path is compared with the expected squared length of a Gaussian vector, that is $n$, since it would be the distribution of the evolution path under random selection. If $\|\pt\|^2$ is greater (respectively lower) than $n$, then the step-size is increased (respectively decreased) following
\begin{equation} \label{eq:sigcsa1}
\st[t+1] = \st \exp \left( \frac{c}{2d_\sigma}\left( \frac{\|\pt[t+1]\|^2}{n} - 1 \right) \right) \enspace ,
\end{equation}
where the damping parameter $d_{\sigma}$ determines how much the step-size can change and can be set here to $d_{\sigma}=1$.

As $[\Ntstar]_i \sim \Nlun$ for $i \geq 3$, we also have $[\pt]_i \sim \Nlun$. It is  convenient in the sequel to also denote by $\pt$ the two dimensional vector $([\pt]_1, [\pt]_2)$. With this (small) abuse of notations, \eqref{eq:sigcsa1} is rewritten as
\begin{equation} \label{eq:sigcsa2}
\st[t+1] = \st \exp \left( \frac{c}{2d_\sigma}\left( \frac{\|\pt[t+1]\|^2 + K_t}{n} - 1 \right) \right) \enspace ,
\end{equation}
with $(K_t)_{t\in\N}$ an i.i.d. sequence of random variables following a chi-squared distribution with $n-2$ degrees of freedom. We shall denote $\stmod$ the multiplicative step-size change $\st[t+1] / \st$, that is the function
\begin{equation} \label{eq:stmod}
\stmod (\pt, \dt, \Wt, K_t) = \exp \left( \frac{c}{2d_\sigma}  \right. \\ \left. \left( \! \frac{\|(1-c)\pt \! +  \sqrt{c(2-c)}\G(\dt, \Wt) \| ^2  +  K_t}{n} \! -\! 1  \right) \! \right) \! \enspace .
\end{equation}

Note that for $c=1$, $\stmod[1]$ is a function of only $\dt$, $\Wt$ and $K_t$ that we will hence denote $\stmod[1](\dt,\Wt,K_t)$.

We prove in the next proposition that for $c<1$ the sequence $(\dt,\pt)_{t \in \N}$ is an homogeneous Markov chain and explicit its update function. In the case where $c=1$ the chain reduces to $\dt$.
  
\begin{proposition} \label{pr:csamarkov}
Consider a $(1,\lambda)$-ES with resampling and cumulative step-size adaptation maximizing the constrained problem \eqref{eq:pbdef}. Take $\dt = g(\Xt)/\st$. The sequence $(\dt, \pt)_{t \in \N}$ is a time-homogeneous Markov chain and
\begin{align} \label{eq:csadelta}
\dt[t+1] &\equald \frac{\dt - \G(\dt, \Wt).\n}{\stmod(\pt, \dt, \Wt, K_t)} \enspace , \\
\ptt &\equald (1-c)\pt + \sqrt{c(2-c)} \G(\dt, \Wt)  \label{eq:csapt} \enspace ,
\end{align}
with $(K_t)_{t \in \N}$ a i.i.d. sequence of  random variables following a chi squared distribution with $n-2$ degrees of freedom, $\G$ defined in Eq.~\eqref{eq:G} \changeaddp{and $\Wt$ defined in Proposition~\ref{pr:g}}.

If $c = 1$ then the sequence $\dmarkov$ is a time-homogeneous Markov chain and
\begin{equation} \label{eq:csadeltanoc}
\dt[t+1] \equald \frac{\dt - \G(\dt,\Wt).\n}{\exp\left(\frac{c}{2d_\sigma}\left(\frac{\|\G(\dt,\Wt)\|^2}{n} - 1 \right)\right)}
\end{equation}
\end{proposition}

\begin{myproof}
With Eq.~\eqref{eq:path} and Eq.~\eqref{eq:NtstarG} we get Eq.~\eqref{eq:csapt}.

From Eq.~\eqref{eq:dt} and Proposition~\ref{pr:g} it follows that
\begin{align*}
\dt[t+1] &= -\frac{\Xt[t+1].\n}{\st[t+1]} \equald -\frac{\Xt.\n + \st \Ntstar.\n}{\st \stmod(\pt, \dt, \Wt, K_t)} \\ &\equald \frac{\dt - \G(\dt, \Wt).\n}{\stmod(\pt, \dt, \Wt, K_t)} \enspace .
\end{align*}
So $(\dt[t+1],\ptt)$ is a function of only $(\dt,\pt)$ and i.i.d. random variables, hence $(\dt,\pt)_{t\in\N}$ is a time-homogeneous Markov chain.

Fixing $c=1$ in \eqref{eq:csadelta} and \eqref{eq:csapt} immediately yields \eqref{eq:csadeltanoc}, and then $\dt[t+1]$ is a function of only $\dt$ and i.i.d. random variables, so in this case $\dmarkov$ is a time-homogeneous Markov chain.
\end{myproof}

As for the constant step-size case, the Markov chain is important when investigating the convergence or divergence of the step size of the algorithm. Indeed from Eq.~\eqref{eq:sigcsa2} we can express $\ln(\st/\st[0])/t$ as
\begin{align} \label{eq:sigrec}
\frac{1}{t}\ln\frac{\st[t]}{\st[0]} &= \frac{c}{2d_\sigma}\left( \frac{ \frac{1}{t} \left(\sum_{i=0}^{t-1}\|\pt[i+1]\|^2 + K_{i}\right)}{n}  - 1 \right)
\end{align}
The right hand side suggests to use the LLN. The convergence of $\ln(\st/\st[0])/t$ to a strictly positive limit (resp. negative) will imply the divergence (resp. convergence) of $\st$ at a geometrical rate.

It turns out that the dynamic of the chain $(\dt, \pt)_{t \in \N}$ looks complex to analyze. Establishing drift conditions looks particularly challenging. We therefore restrict the rest of the study to the more simple case where $c=1$, hence the Markov chain of interest is $\dmarkov$. Then \eqref{eq:sigrec} becomes
\begin{equation} \label{eq:sigrec1}
\frac{1}{t}\ln\frac{\st[t]}{\st[0]} \equald \frac{c}{2d_\sigma}  \left( \frac{ \frac{1}{t}\sum_{i=0}^{t-1}  \|\G(\dt[i],\Wt[i])\|^2  +  K_i}{n}  - 1 \right) \enspace .
\end{equation}

To apply the LLN we will need the Markov chain to be Harris positive, and the properties mentioned in the following lemma. 


\begin{lemma}[{\citealt[Proposition~7]{CA2014irrapfel}}] \label{lm:csaproperties}
Consider a $(1,\lambda)$-ES with resampling and cumulative step-size adaptation maximizing the constrained problem \eqref{eq:pbdef}. For $c=1$ the Markov chain $\dmarkov$ from Proposition~\ref{pr:csamarkov} is $\psi$-irreducible, strongly aperiodic, and compact sets \adp{of $\R^*_+$} are small sets for this chain.
\end{lemma}

\newcommand{\h}{\bs{h}}
\renewcommand{\u}{\bs{u}}

\anne{I know that I suggested to put the proof within the latex, but now I see that it is not really understandable and that anyhow we refer to lemmas and definitions from the Bernoulli to be sumitted paper. Hence, given that the arxiv paper should stay forever on the web, I would say it's fine to only refer to the arxiv paper (of the first version of the paper to be submitted to bernoulli) and skip the proof.}
\mathnote{
\begin{proof}
To prove the lemma we use \cite[Theorem~1]{CA2014irrapfel}. We first check that the model of this  theorem is respected.

Let $F : \R\times \R^2 \to \R$ be the function
\begin{equation}
F(\delta, \w) = \frac{\delta - \w.\n}{\exp\left(\frac{c}{2d_\sigma}\left(\frac{\|\w\|^2}{n} - 1 \right)\right)} \enspace .
\end{equation}
Then according to Proposition~\ref{pr:csamarkov}, $\dt[t+1] \equald F(\dt, \G(\dt, \Wt))$, where $(\Wt)_{t \in \N}$ is a i.i.d. sequence of random vectors defined in Proposition~\ref{pr:g}. The function $F$ is $C^1$, and the random variable $\G(\delta, \Wt)$ admits a density $\pstar$ due to Lemma~\ref{lm:pstar} and Proposition~\ref{pr:g} which is lower semi-continuous.

We define $F^1$ as $F$, and for $t \in \N^*$ we inductively define
\begin{equation}
F^{t+1}(\delta, \w_1, \ldots, \w_{t+1}) := F^t( F(\delta, \w_1), \w_2, \ldots, \w_{t+1}) \enspace .
\end{equation}
Similarly, $p_{\delta, 1}^\star := \pstar$ and inductively
\begin{equation}
p_{\delta, t+1}^\star( \w_1, \ldots, \w_{t+1}) := p_{\delta, t}^\star(\delta, \w_1, \ldots, \w_t) \pstar[F^t(\delta, \w_1, \ldots, \w_t)](\w_{t+1}) \enspace,
\end{equation}
and for $\delta \in \R$ we define $O_{\delta,t}$ as
\begin{equation}
O_{\delta,t} := \lbrace (\w_1, \ldots \w_t) \in \R^{2\times t} | p_{\delta, t}^\star(\delta, \w_1, \ldots, \w_t) > 0 \rbrace  \enspace ,
\end{equation}
and $O_\delta$ as $O_{\delta,1} = \lbrace \w \in \R^{2} | \delta - \w.\n > 0\rbrace$.

\newcommand{\0}{\bs{0}}

We now prove that the point $\delta^* \in \R_+^*$ is a strongly globally attracting state, i.e. for all $\delta_0 \in \R_+^*$ and $\epsilon \in \R_+^*$ there exists $t_0 \in \N^*$ such that for all $t \geq t_0$ there exists $\w \in O_{\delta_0, t}$ such that $F^t(\delta_0, \w) \in B(\delta^*, \epsilon)$. Let $\delta_0 \in \R_+^*$. Let $k \in \N^*$ be such that $\delta_0 \exp(ck/(2d_\sigma)) > \delta^*$. We take $\w_i =  \0$ for all $i \in [1..k]$ and define $\delta_k := F^k(\delta_0, \w_1, \ldots, \w_k)$. By construction of $k$, we have $\delta_k > \delta^*$. Now, take $\u = (-1, -1)$ and note that the limit $\lim_{\alpha \to +\infty} F(\delta_k, \alpha\u) = 0$. Since the function $F$ is continuous and that $F(\delta_k, \0) > \delta_k$, this means that the function $\alpha \mapsto F(\delta_k, \alpha \u)$ is surjective in $(0, \delta_k)$. And since $\delta^* < \delta_k$, there exists $\alpha_0$ such that $F(\delta_k, \alpha_0 \u) = \delta^*$. Now let $\w$ denote $(\w_1, \ldots, \w_k, \alpha_0 \u)$, and note that $\alpha \u \in O_{\delta} = \lbrace \bs{v} \in \R^2 | \delta - [\bs{v}]_1 \cos \theta - [\bs{v}]_2 \sin \theta > 0 \rbrace $ for all $\alpha \in \R_+$ and all $\delta \in \R_+^*$; hence $\alpha_0 \u \in O_{\delta_k}$ and $\w_i = 0\u \in O_{\delta}$ for all $\delta \in \R_+^*$. Therefore $\w \in O_{\delta_0, k+1}$ and $F^{k+1}(\delta_0, \w) = \delta^*$. As the proof stands for all $k$ large enough, $\delta^*$ is a strongly globally attractive state.

We will now show that $F(\0, \cdot)$ is a submersion at some point $\w \in \R^2$. To do so we calculate the differential $D_\w F(0, \cdot)$ of $F(0, \cdot)$ at $\w$:
\begin{align*}
F(0, \w + \h) &= g(\w + \h) \exp \left( -\frac{1}{2d_\sigma} \left( \frac{\| \w + \h \|^2}{ n } - 1\right) \right) \\
 &= g(\w + \h) \exp \left( -\frac{1}{2d_\sigma} \left( \frac{\| \w \|^2 + 2\w.\h  + \| \h \|^2}{ n } - 1\right) \right) \\
 &= g(\w + \h) \exp\left( -\frac{1}{2d_\sigma} \left( \frac{\| \w \|^2 }{ n } - 1\right) \right) \exp\left(- \frac{1}{2 d_\sigma} \left( \frac{2\w.\h + \|\h\|^2}{n}  \right) \right) \\
 &= g(\w + \h) \exp\left( -\frac{1}{2d_\sigma} \left( \frac{\| \w \|^2 }{ n } - 1\right) \right)\left(1 - \frac{2 \w.\h}{2 d_\sigma n}  + o(\|\h\|)\right) \\
 &= F(0, \w) - F(0, \w)\frac{\w.\h}{d_\sigma n} +  g(\h) \exp \left(-\frac{1}{2d_\sigma} \left( \frac{\|\w\|^2}{n} - 1 \right) \right) + o(\|\h\|) \enspace .
\end{align*}
Hence for $\w = (-\sqrt{n}, 0, \ldots, 0)$ and $\h = (0, \alpha, 0, \ldots, 0)$, $D_\w F(0, \cdot) (\h) = -\alpha \sin \theta \exp(0)$. Hence for $\alpha$ spanning $\R$, $D_\w F(0, \cdot) (\h)$ spans $\R$ such that the image of $D_\w F(0, \cdot)$ equals $\R$, i.e. $D_\w F(0, \cdot)$ is surjective meaning $F(0, \cdot)$ is a submersion at $\w$. According to \cite[Lemma~5]{CA2014irrapfel} this means there exists $N$ an open neighbourhood of $(0, \w)$ such that for all $(\delta, \u) \in N$, $F(\delta, \cdot)$ is a submersion at $\u$. So for $\delta^* \in \R_+^*$ small enough, $F(\delta^*, \cdot)$ is a submersion at $\w \in O_{\delta^*}$.

Adding this with the fact that $\delta^*$ is a strongly globally attracting state, we can then apply \cite[Theorem~1]{CA2014irrapfel} which concludes the proof.

\end{proof}
}

\newcommand{\wt}[1][t]{\bs{w}_{#1}}
We believe that the latter result can be generalized to the case $c<1$ if for any $(\delta_0, \p_0) \in \ddomain\times \R^n$ there exists $t_{\delta_0, \p_0}$ such that for all $t \geq t_{\delta_0, \p_0}$ there exists a path of events of length $t$ from $(\delta_0, \p_0)$ to \rep{any point of} the set $\rep{[}\adp{(}0,M] \times B(\bs{0},r)$ \adp{for $M> 0$ and $r>0$ small enough}.

To show the Harris positivity of $\dmarkov$ we \adp{will use the drift function $V: \delta \in \R_+^* \mapsto \delta^\alpha + \delta^{-\alpha}$. From the definition of the drift operator $\Delta V$ in \eqref{eq:driftop} and the update of $\dt$ in \eqref{eq:csadeltanoc}, we then have
\begin{equation} \label{eq:driftopcumul}
\Delta V(\delta) = \E\left( \frac{\left( \delta - \G(\delta, \W).\n\right)^\alpha}{\stmod[1](\delta,\W,K)^\alpha}\right) + \E\left( \frac{\stmod[1](\delta,\W,K)^\alpha}{\left( \delta - \G(\delta, \W).\n\right)^\alpha}\right) - \delta^\alpha - \delta^{-\alpha} \enspace .
\end{equation}
  To verify the drift condition of \eqref{eq:V4}, using the fact from Lemma~\ref{lm:csaproperties} that for $0<m<M$ the compact $[m,M]$ is a small set, it is sufficient to show that the limits of $\Delta V / V(\delta)$ in $0$ and $\infty$ is negative. These limits \adp{will result from the limits studied}\rep{ are studied} in the following lemma corresponding the the decomposition in \eqref{eq:driftopcumul}. } \rep{ first need to study the behaviour of the drift operator \rep{ we want to use}\del{$\Delta V$ defined in \eqref{eq:driftop}} when $\delta \rightarrow +\infty$, that is far from the constraint. Then, intuitively, as $[\Ntstar]_2$ would not be influenced by the resampling anymore, it would be distributed as a random normal variable, and $[\Ntstar]_1$ would be distributed as the last order statistic of $\lambda$ normal random variables. \adp{We also need to control the drift operator when $\delta \to 0$. The limits when $\delta \to + \infty$ and $\delta \to 0$ of the drift operator  are studied }\rep{ This is used} in the following technical lemma.}

\begin{lemma} \label{lm:farcsa}
For $\alpha > 0$ \adp{small enough}
\begin{align} \label{eq:farcsa}
\adp{\frac{1}{\delta^\alpha + \delta^{-\alpha} }} \E\left(\frac{\left(\delta -  \G(\delta, \W).\n \right)^\alpha} { \rep{\delta^\alpha} \stmod[1](\delta,\W, K)^{\adp{\alpha}}}\right) &\underset{\delta \rightarrow +\infty}{\longrightarrow} E_1 E_2 E_3 \adp{< \infty} \\
\adp{\frac{1}{\delta^\alpha + \delta^{-\alpha}} \E\left(\frac{\left(\delta -  \G(\delta, \W).\n \right)^\alpha} { \stmod[1](\delta,\W, K)^\alpha}\right)} &\adp{\underset{\delta \rightarrow 0}{\longrightarrow} 0 }\\
\adp{\frac{1}{\delta^\alpha + \delta^{-\alpha}} \E\left(\frac { \stmod[1](\delta,\W, K)^\alpha}{\left(\delta -  \G(\delta, \W).\n \right)^\alpha} \right)} &\adp{\underset{\delta \rightarrow +\infty}{\longrightarrow} 0 }\\
\adp{\frac{1}{\delta^\alpha + \delta^{-\alpha}} \E\left(\frac { \stmod[1](\delta,\W, K)^\alpha}{\left(\delta -  \G(\delta, \W).\n \right)^\alpha}\right)}&\adp{\underset{\delta \rightarrow 0}{\longrightarrow} 0} \enspace ,
\end{align} 
where $E_1 = \E ( \exp (-\frac{\alpha}{2d_\sigma n}(\Nlambda^2 - 1 )))$,
$E_2 = \E ( \exp (-\frac{\alpha}{2d_\sigma n}(\Nlun^2 - 1 )))$, and $E_3 = \E ( \exp (-\frac{\alpha}{2d_\sigma n}(K - (n-2) )))$; where $\G$ is the function defined in Eq.~\eqref{eq:G} and $\stmod[1]$ is defined in Eq.~\eqref{eq:stmod} (for $c=1$), $K$ is a random variable following a chi-squared distribution with $n-2$ degrees of freedom and $\W \sim (\Ud_{[0,1]}, \Nlun)^\lambda$ is a random vector.
\end{lemma}

The proof of this lemma consists in applications of Lebesgue's dominated convergence theorem, and can be found in the appendix.

We now prove the Harris positivity of $\dmarkov$ by proving a stronger property, namely the geometric ergodicity that we show using the drift inequality \eqref{eq:V4}.

\begin{proposition} \label{pr:csaerg}
Consider a $(1,\lambda)$-ES with resampling and cumulative step-size adaptation maximizing the constrained problem \eqref{eq:pbdef}. For $c=1$ the Markov chain $\dmarkov$ from Proposition~\ref{pr:csamarkov} is $V$-geometrically ergodic with $V : \delta \in \ddomain \mapsto \delta^\alpha + \rep{1}\adp{\delta^{-\alpha}}$ for $\alpha \adp{>0}$ small enough, and  positive Harris with invariant measure $\pi_1$.
\end{proposition}

\begin{myproof}
Take $V$ the positive function $V(\delta) =  \delta^\alpha + \rep{1}\adp{\delta^{-\alpha}}$ (the parameter $\alpha$ is strictly positive  and will be specified later), $\W \sim (\Ud_{[0,1]}, \Nlun)^\lambda$ a random vector and $K$ a random variable following a chi squared distribution with $n-2$ degrees of freedom. \adp{We first study $\Delta V/V(\delta)$ when $\delta \to + \infty$.} From Eq.~\adp{\eqref{eq:driftopcumul}}\rep{\eqref{eq:csadeltanoc}} we then have the following drift quotient 
\begin{align} \label{eq:csadeltavv}
\frac{\Delta V(\delta)}{V(\delta)\re{-1}} &=  \adp{\frac{1}{V(\delta)}}\E \left( \frac{( \delta - \G(\delta,\W).\n )^\alpha}{\rep{\delta^\alpha} \stmod[1](\delta,\W, K)^{\adp{\alpha}}}\right)   
\adp{ +  \adp{\frac{1}{V(\delta)}}\E\left( \frac{\stmod[1](\delta,\W, K)^\alpha}{(\delta - \G(\delta,\W).\n )^\alpha} \right)} - 1  \enspace ,
\end{align}
with $\stmod[1]$ defined in Eq.~\eqref{eq:stmod} and $\G$ in Eq.~\eqref{eq:G}.  From Lemma~\ref{lm:farcsa}, following the same notations than in the lemma, when $\delta \rightarrow +\infty$ \adp{and if $\alpha>0$ is small enough,} the right hand side of the previous equation converges to $E_1 E_2 E_3-1$. With Taylor series \rep{$E_1 = \E( \sum_{k \in \N} (-\alpha/(2d_\sigma n)(\Nlambda^2 - 1))^k / k!)$. As $\E( \sum_{k \in \N} |-\alpha/(2d_\sigma n)(\Nlambda^2 - 1))^k / k!|) = \int \exp | \alpha/(2 d_\sigma n)(x^2 - 1)| \lambda \varphi(x) \Phi(x)^{\lambda -1} \mathrm{d}x$ and that}
\adp{
\begin{align*}
E_1 &= \E \left( \sum_{k \in \N} \frac{\left(-\frac{\alpha}{2d_\sigma n}\left(\Nlambda^2 - 1\right)\right)^k }{k!}\right)  \enspace . 
\end{align*}
Furthermore, as the density of $\Nlambda$ at $x$ equals to $\lambda \varphi(x) \Phi(x)^{\lambda-1}$ and that}
 $\exp | \alpha/(2 d_\sigma n)(x^2 - 1)| \lambda \varphi(x) \Phi(x)^{\lambda -1} \leq \lambda \adp{\exp(1/(2 d_\sigma n))} \exp( \alpha/(2 d_\sigma n) x^2 - x^2/2)$  which for $\alpha$ small enough is integrable, \adp{hence }
\adp{
\begin{equation*}
\E \left( \sum_{k \in \N} \frac{\left|-\frac{\alpha}{2d_\sigma n}\left(\Nlambda^2 - 1\right)\right|^k }{k!}\right) = \int_\R \exp \left| \frac{\alpha}{2 d_\sigma n}\left(x^2 - 1\right)\right| \lambda \varphi(x) \Phi(x)^{\lambda -1} \mathrm{d}x < \infty \enspace .
\end{equation*}
Therefore
}
we can use Fubini's theorem to invert series \adp{(which are integrals for the counting measure)}  and integral. The same reasoning holding for $E_2$ and $E_3$ (for $E_3$ with the chi-squared distribution we need $\alpha/(2d_\sigma n)x - x/2 \leq 0$ \adp{for all $x\geq 0$}) we have
\begin{align*}
\lim_{\delta \rightarrow + \infty} \frac{\Delta V}{V\rep{-1}}(\delta) =  \left( \! 1 \! - \! \frac{\alpha}{2d_\sigma n}\! \E(\Nlambda^2 \!\! - \! 1) \! + \! o(\alpha) \! \right) \!\! \left(\! 1 \! - \!\! \frac{\alpha}{2d_\sigma n}\E(\Nlun^2 \!\! - \! 1) \! + \! o(\alpha)\! \! \right) \\ \left( 1 - \frac{\alpha}{2d_\sigma n}\E(\chisq - (n-2)) + o(\alpha)  \right) - 1 \enspace ,
\end{align*}
and as $\E(\Nlun^2) = 1$ and $\E(\chisq) = n-2$
\begin{equation*}
\rep{\lim_{\delta \rightarrow + \infty} \frac{\Delta V}{V-1}(\delta) =} \lim_{\delta \rightarrow + \infty} \frac{\Delta V}{V}(\delta) = - \frac{\alpha}{2d_\sigma n}\E\left(\Nlambda^2 - 1\right) + o(\alpha) \enspace .
\end{equation*}
From \cite{cahPPSN12} if $\lambda > 2$ then $\E(\Nlambda^2) > 1$. Therefore, for $\alpha$ small enough, we have $\lim_{\delta \rightarrow +\infty} \frac{\Delta V}{V}(\delta) < 0$  so there exists $\epsilon_{\adp{1}} > 0$\rep{, $b \in \R$} and $M > 0$ such that  \adp{$\Delta V (\delta) \leq -\epsilon_{\adp{1}} V(\delta)$ whenever $\delta > M$.}

\adp{ Similarly, when $\alpha$ is small enough, using Lemma~\ref{lm:farcsa}, $\lim_{\delta \to 0} \E((\delta-\G(\delta, \W))^\alpha/\stmod[1](\delta,\W,K)^\alpha)/V(\delta) = 0$ and $\lim_{\delta \to 0} \E(\stmod[1](\delta,\W,K)^\alpha/(\delta-\G(\delta, \W))^\alpha)/V(\delta) = 0$. Hence using \eqref{eq:csadeltavv},  $\lim_{\delta \to 0} \Delta V (\delta) / V (\delta) = -1$. So there exists $\epsilon_2$ and $m > 0$ such that $\Delta V(\delta) \leq -\epsilon_2 V(\delta)$ for all $\delta \in (0,m)$. And since $\Delta V(\delta)$ and $V(\delta)$ are bounded functions on compacts of $\R_+^*$, there exists $b \in \R$ such that
}
\begin{equation*}
\Delta V(\delta) \leq -\rep{\epsilon}\adp{\min(\epsilon_1, \epsilon_2)} V(\delta) + b \1_{[\rep{0}\adp{m},M]}(\delta) \enspace .
\end{equation*}

With Lemma~\ref{lm:csaproperties}, $[\rep{0}\adp{m},M]$ is a small set, and $\dmarkov$ is a $\psi$-irreducible aperiodic Markov chain. So $\dmarkov$ satisfies the assumptions of \cite[Theorem 15.0.1]{markovtheory}, which proves the proposition.

\end{myproof}

The same results for $c<1$ are difficult to obtain, as then both $\dt$ and $\pt$ must be controlled together. For $\pt = 0$ and $\dt \geq M$, $\|\ptt\|$ and $\dt[t+1]$ will in average increase, so either we need that $[M, +\infty)\times B(\bs{0},r)$ is a small set (although it is not compact), or we need to look $\tau$ steps in the future with $\tau$ large enough to see $\dt[t+\tau]$ decrease for all possible values of $\pt$ outside of a small set.

Note that although in Proposition~\ref{pr:csterg} and Proposition~\ref{pr:csaerg} we show the existence of a stationary measure for $\dmarkov$, these are not the same measures, and not the same Markov chains as they have different update rules (compare Eq.~\eqref{eq:cstmarkov} and Eq.~\eqref{eq:csadelta}).
The chain $\dmarkov$ being Harris positive we may now apply a LLN to Eq.~\eqref{eq:sigrec1} to get an exact expression of the divergence/convergence rate of the step-size.

\begin{theorem} \label{th:csa}
Consider a $(1,\lambda)$-ES with resampling and cumulative step-size adaptation maximizing the constrained problem \eqref{eq:pbdef}, and for $c=1$ take $\dmarkov$ the Markov chain from Proposition~\ref{pr:csamarkov}.
Then the step-size diverges or converges geometrically in probability
\begin{equation} \label{eq:csa}
\frac{1}{t} \ln \left( \frac{\st}{\st[0]} \right) \overset{P}{\underset{t \rightarrow \infty}{\longrightarrow}}  \frac{1}{2d_\sigma n}\left(  \E_{\pi_1\otimes\mu_\W} \left(\|\G\left(\delta, \W\right)\|^2\right) - 2 \right) \enspace ,
\end{equation}
and in expectation
\begin{equation} \label{eq:csaesp}
\E\left(\ln\left(\frac{\st[t+1]}{\st}\right)\right) \underset{t \rightarrow +\infty}{\longrightarrow} \frac{1}{2d_\sigma n}\left(  \E_{\pi_1\otimes\mu_\W} \left(\|\G\left(\delta, \W\right)\|^2\right) - 2 \right)
\end{equation}
with $\G$ defined in \eqref{eq:G} and $\W = (\bs{W}^i)_{i \in [1..\lambda]}$ where $(\bs{W}^i)_{i \in [1..\lambda]}$ is an i.i.d. sequence such that $\bs{W}^i \sim (\Ud_{[0,1]}, \Nlun)$, $\mu_\W$ is the probability measure of $\W$ and $\pi_1$ is the invariant measure of $\dmarkov$ whose existence is proved in Proposition~\ref{pr:csaerg}.

Furthermore, the change in fitness value $f(\Xtt) - f(\Xt)$ diverges or converges geometrically in probability
\begin{equation} \label{eq:csafx}
\frac{1}{t} \ln \left| \frac{  f(\Xtt) - f(\Xt) }{\st[0]}  \right|  \overset{P}{\underset{t \rightarrow \infty}{\longrightarrow}}  \frac{1}{2d_\sigma n}\left(  \E_{\pi_1\otimes\mu_\W} \left(\|\G\left(\delta, \W\right)\|^2\right) - 2 \right) \enspace .
\end{equation}

\end{theorem}

\begin{myproof}

From Proposition~\ref{pr:csaerg} the Markov chain $\dmarkov$ is Harris positive, and since $(\Wt)_{t \in \N}$ is i.i.d., the chain $(\dt, \Wt)_{t \in \N}$ is also Harris positive with invariant probability measure $\pi_1 \times \mu_\W$, so to apply the Law of Large Numbers of \cite[Theorem 17.0.1]{markovtheory} to Eq.~\eqref{eq:sigrec} we only need the function $(\delta,\bs{w}) \mapsto \|G(\delta,\bs{w})\|^2 + K$ to be $\pi_1\times\mu_\W$-integrable.

Since $K$ has chi-squared distribution with $n-2$ degrees of freedom, $\E_{\pi_1\times\mu_\W}(\|G(\delta,\W)\|^2 + K)$ equals to $\E_{\pi_1\times\mu_\W}(\|G(\delta,\W)\|^2) + n-2$. With Fubini-Tonelli's theorem, $\E_{\pi_1\times\mu_\W}(\|\G(\delta,\W)\|^2)$ is equal to  $\E_{\pi_1}(\E_{\mu_\W}(\|\G(\delta,\W)\|^2 ))$. From Eq.~\eqref{eq:pstar} and from the proof of Lemma~\ref{lm:farNtstar} the function $\x \mapsto \|\x\|^2\pstar(\x)$ converges simply to $\|\x\|^2 {p_{\Nlambda}} ([\x]_1) \varphi([x]_2)$ while being dominated by $\lambda / \Phi(0) \exp(-\|\x\|^2)$ which is integrable. Hence we may apply Lebesgue's dominated convergence theorem showing that the function $\delta \mapsto \E_{\mu_{\W}}(\|\G(\delta,\W)\|^2)$ is continuous and has a finite limit and is therefore bounded by a constant $M_{\G^2}$. As the measure $\pi_1$ is a probability measure (so $\pi_1(\R) = 1$), $\E_{\pi_1}(\E_{\mu_\W}(\|\G(\delta,\W)\|^2 | \dt = \delta)) \leq M_{\G^2} < \infty$. Hence we may apply the Law of Large Numbers
\begin{equation*}
\sum_{i=0}^{t-1} \! \frac{\|\G(\dt[i], \Wt[i])\|^2 \! + \! K_i}{t} \! \overset{a.s}{\underset{t \rightarrow \infty}{\longrightarrow}} \! \E_{\pi_1 \times \mu_\W} \! \left( \|\G(\delta,\W)\|^2 \right) + n-2 \enspace .
\end{equation*}
Combining this equation with Eq.~\eqref{eq:sigrec1} yields Eq.~\eqref{eq:csa}.

From Proposition~\ref{pr:g}, \eqref{eq:path} for $c=1$ and \eqref{eq:sigcsa2}, $\ln(\stt/\st) \equald 1/(2d_\sigma n) (\|\G(\dt,\Wt)\|^2 + \chisq - n)$ so $\E(\ln(\stt/\st) | (\dt[0], \st[0]) ) =  1/(2d_\sigma n) ( \E(\|\G(\dt,\Wt)\|^2 | (\dt[0], \st[0])) - 2)$. As $\|\G\|^2$ is integrable with Fubini's theorem $\E(\|\G(\dt, \Wt)\|^2 | (\dt[0], \st[0]) ) = \int_{\ddomain}  \E_{\mu_\W}(\|\G( \y, \W)\|^2) P^t(\dt[0], d\y) $, so $ \E(\|\G(\dt, \Wt)\|^2 | (\dt[0], \st[0])) - \E_{\pi_1 \times \mu_\W}( \|\G(\delta, \W) \|^2) = \int_{\ddomain} \E_{\mu_\W}( \|\G( \y, \W)\|^2 ) (P^t(\x/\sigma, d\y) - \pi_1(d\y))  $. According to Proposition~\ref{pr:csaerg} $\dmarkov$ is $V$-geometrically ergodic with $V : \delta \mapsto \delta^\alpha + \rep{1/}\adp{\delta^{-\alpha}}$, so there exists $M_\delta$ and $r > 1$ such that $\|P^t(\delta, \cdot) - \pi_1 \|_V \leq M_\delta r^{-t} $. We showed that the function $\delta \mapsto \E(\|\G(\delta, \W)\|^2)$ is bounded, so since $V(\delta) \geq 1$ for all $\delta \in \ddomain$  there exists $k$ such that $\E_{\mu_\W}(\|\G(\delta, \W)\|^2) \leq k V(\delta)$ for all $\delta$. Hence $| \int \E_{\mu_\W}(\|\G( x, \W)\|^2) (P^t(\delta, dx) - \pi_1(dx)) | \leq  k\|P^t(\delta, \cdot) - \pi_1 \|_V \leq k M_\delta r^{-t} $. And therefore $|\E(\|\G(\dt, \Wt)\|^2 | (\dt[0], \st[0])) - \E_{\pi_1 \times \mu_\W}(\|\G(\delta, \W))\|^2)  \leq k M_\delta r^{-t}$ which converges to $0$ when $t$ goes to infinity, which shows Eq.~\eqref{eq:csaesp}.


For \eqref{eq:csafx} we have that $\Xtt - \Xt \equald \st \G(\dt, \Wt)$ so $\changeaddp{(}1/t\changeaddp{)}\ln|(f(\Xtt)-f(\Xt))/\st[0]| \equald \changeaddp{(}1/t\changeaddp{)}\ln(\st/\st[0]) + \changeaddp{(}1/t\changeaddp{)}\ln|f(\G(\dt, \Wt))/\st[0]| $. From \eqref{eq:pstar1}, since $1/2 \leq \Phi(x) \leq 1$ for all $x \geq 0$ and that $\cdfnone(x) \leq 1$,  the probability density function of $f(\G(\dt,\Wt)) = [\G(\dt,\Wt)]_1$ is dominated by $2 \lambda \varphi(x)$. Hence 
\begin{align*}
\Pr(\ln|[\G(\delta, \W)]_1|/t \geq \epsilon) &\leq \int_{\R} \1_{[\epsilon t, +\infty)}(\ln|x|)  2\lambda \varphi(x) dx \\
&\leq \int_{\exp(\epsilon t)}^{+\infty} 2\lambda \varphi(x)dx + \int_{-\infty}^{-\exp(\epsilon t)} 2\lambda \varphi(x)dx
\end{align*}
For all $\epsilon > 0$ since $\varphi$ is integrable with the dominated convergence theorem both members of the previous inequation converges to $0$ when $t \rightarrow \infty$, which shows that $\ln|f(\G(\dt,\Wt))|/t$ converges in probability to $0$. Since $\ln(\st/\st[0])/t$ converges in probability to the right hand side of \eqref{eq:csafx} we get \eqref{eq:csafx}.
\end{myproof}

If, for $c<1$, the chain $(\dt, \pt)_{t \in \N}$ was positive Harris with invariant measure $\pi_c$ and $V$-ergodic such that $\|\ptt\|^2$ is dominated by $V$ then we would obtain similar results with a convergence/divergence rate equal to $c /(2 d_\sigma n)(  \E_{{\pi_c}\otimes\mu_\W} (\| \bs{p} \|^2 ) - 2 )$.

If the sign of the RHS of Eq.~\eqref{eq:csa} is strictly positive then the step size diverges geometrically. The Law of Large Numbers entails that Monte Carlo simulations will converge to the RHS of Eq.~\ref{eq:csa}, and the fact that the chain is $V$-geometrically ergodic (see Proposition~\ref{pr:csaerg}) means sampling from the $t$-steps transition kernel $P^t$ will get close exponentially fast to sampling directly from the stationary distribution $\pi_1$. We could apply a Central Limit Theorem for Markov chains~\cite[Theorem 17.0.1]{markovtheory}, and get an approximate confidence interval for $\ln(\st/\st[0])/t$, given that we find a function $V$ for which the chain $(\dt, \Wt)_{t \in \N}$ is $V$-uniformly ergodic and such that $\|\G(\delta,\w)\|^4 \leq V(\delta, \w)$. The question of the sign of $\lim_{t \to +\infty} f(\Xt)-f(\Xt[0])$ is not adressed in Theorem~\ref{th:csa}, but simulations indicate that for $d_\sigma \geq 1$ the probability that $f(\Xt)>f(\Xt[0])$ converges to $1$ as $t \to +\infty$. For low enough values of $d_\sigma$ \adp{and of $\theta$} this probability appears to converge to $0$.
\alex{do I add plots to support this ? In this case, what do I plot ? I think it is interesting to see the probability against time. The speed of convergence to $1$ depends on the parameters.}\niko{I would say only after you finished your thesis. }

\begin{figure}
\begin{center}
\includegraphics[width=0.70\textwidth]{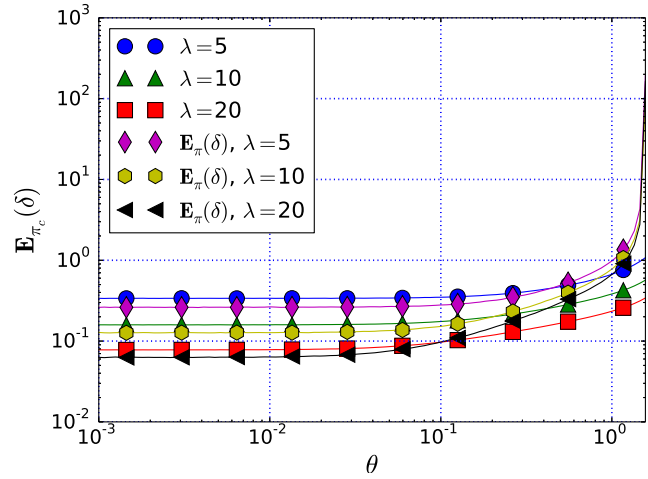}
\end{center}
\caption{Average normalized distance $\delta$ from the constraint for the $(1,\lambda)$-CSA-ES plotted against the constraint angle $\theta$, for $\lambda \in \lbrace 5, 10, 20 \rbrace$, $c=1/\sqrt{2} $\changeaddp{, $d_\sigma = 1$} and dimension $2$.}
\label{fg:csadvst}
\end{figure}

\begin{figure} 
\begin{center}
\includegraphics[width=0.70\textwidth]{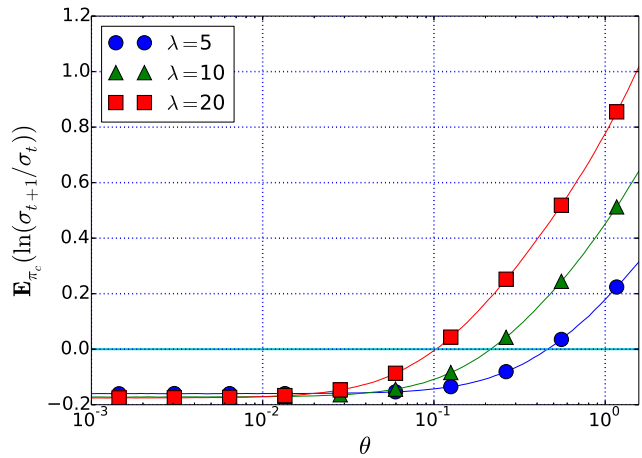}
\end{center}
\caption{Average of the logarithmic adaptation response $\Delta_t = \ln(\sigma_{t+1}/\sigma_t)$ for the $(1,\lambda)$-CSA-ES plotted against the constraint angle $\theta$, for $\lambda \in \lbrace 5, 10, 20 \rbrace$, $c = 1/\sqrt{2}$\changeaddp{, $d_\sigma = 1$} and dimension $2$. Values below zero (straight line) indicate premature convergence.
}
\label{fg:csasvst}
\end{figure}

As in Fig.~\ref{fg:cstdvst} we simulate the Markov chain $(\dt, \pt)_{t \in \N}$ defined in Eq.~\eqref{eq:csadelta} to obtain Fig.~\ref{fg:csadvst} after an average of $\dt$ over $10^6$ time steps. \adp{Assuming that the Markov chain $\dpmarkov$ admits an invariant probability measure $\pi_c$, t}he expected value $\E_{\pi_{\rep{1}\adp{c}}}(\delta)$ shows the same dependency in $\lambda$ \adp{as}\rep{that} in the constant case. With \changeremp{higher}\changeaddp{larger} population size, the algorithm follows the constraint from closer, as better samples are available closer to the constraint, which a \changeremp{higher}\changeaddp{larger} population helps to find. The difference between $\E_{\pi_c}(\delta)$ and $\E_\pi(\delta)$ appears small except for large values of the constraint angle. When $\E_\pi(\delta) > \E_{\pi_c}(\delta)$ we observe on Fig.~\ref{fg:csasvst} that $\E_{\pi_c}(\ln(\stt/\st)) > 0$.

\begin{figure} 
\begin{center}
\includegraphics[width=0.70\textwidth]{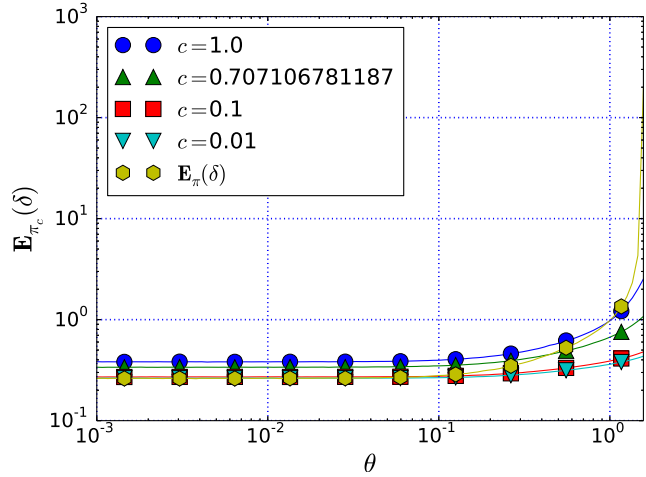}
\end{center}
\caption{Average normalized distance $\delta$ from the constraint for the $(1,\lambda)$-CSA-ES plotted against the constraint angle $\theta$ with $c \in \lbrace 1, 1/\sqrt{2}, 0.1, 0.01 \rbrace$ and for constant step-size, where $\lambda = 5$\changeaddp{, $d_\sigma = 1$} and dimension $2$.}
\label{fg:csadvstcs}
\end{figure}

In Fig.~\ref{fg:csadvstcs} the average of $\dt$ over $10^6$ time steps is again plotted with $\lambda = 5$, this time for different values of the cumulation parameter, and compared with the constant step-size case. A lower value of $c$ makes the algorithm follow the constraint from closer. When $\theta$ goes to $0$ the value $\E_{\pi_c}(\delta)$ converges to a constant, and $\lim_{\theta \rightarrow 0}\E_{\pi}(\delta)$ for constant step-size seem to be $\lim_{\theta \rightarrow 0}\E_{\pi_c}(\delta)$ when $c$ goes to $0$. As in Fig.~\ref{fg:csadvst} the difference between $\E_{\pi_c}(\delta)$ and $\E_\pi(\delta)$ appears small except for large values of the constraint angle. This suggests that the difference between the distributions $\pi$ and $\pi_c$ is small. Therefore the approximation made in \cite{arnold2011behaviour} where $\pi$ is used instead of $\pi_c$ to estimate $\ln(\stt/\st)$ is accurate for not too large values of the constraint angle.

In Fig.~\ref{fg:csasvst}, corresponding to the LHS of Eq.~\eqref{eq:csa}, \adp{the adaptation response $\Delta_t := \ln(\sigma_{t+1}/\sigma_t)$ is averaged}\del{adaptation response $\Delta_t$ is averaged}\rep{left hand side of Eq.~\eqref{eq:csa} is simulated} over $10^6$ time steps \adp{and plotted} against the constraint angle $\theta$ for different population sizes. If the value is below zero the step-size converges, which means a premature convergence of the algorithm. We see that a \changeremp{higher}\changeaddp{larger} population size helps to achieve a faster divergence rate and for the step-size adaptation to succeed for a wider interval of values of $\theta$. 

\begin{figure} 
\begin{center}
\includegraphics[width=0.49\textwidth]{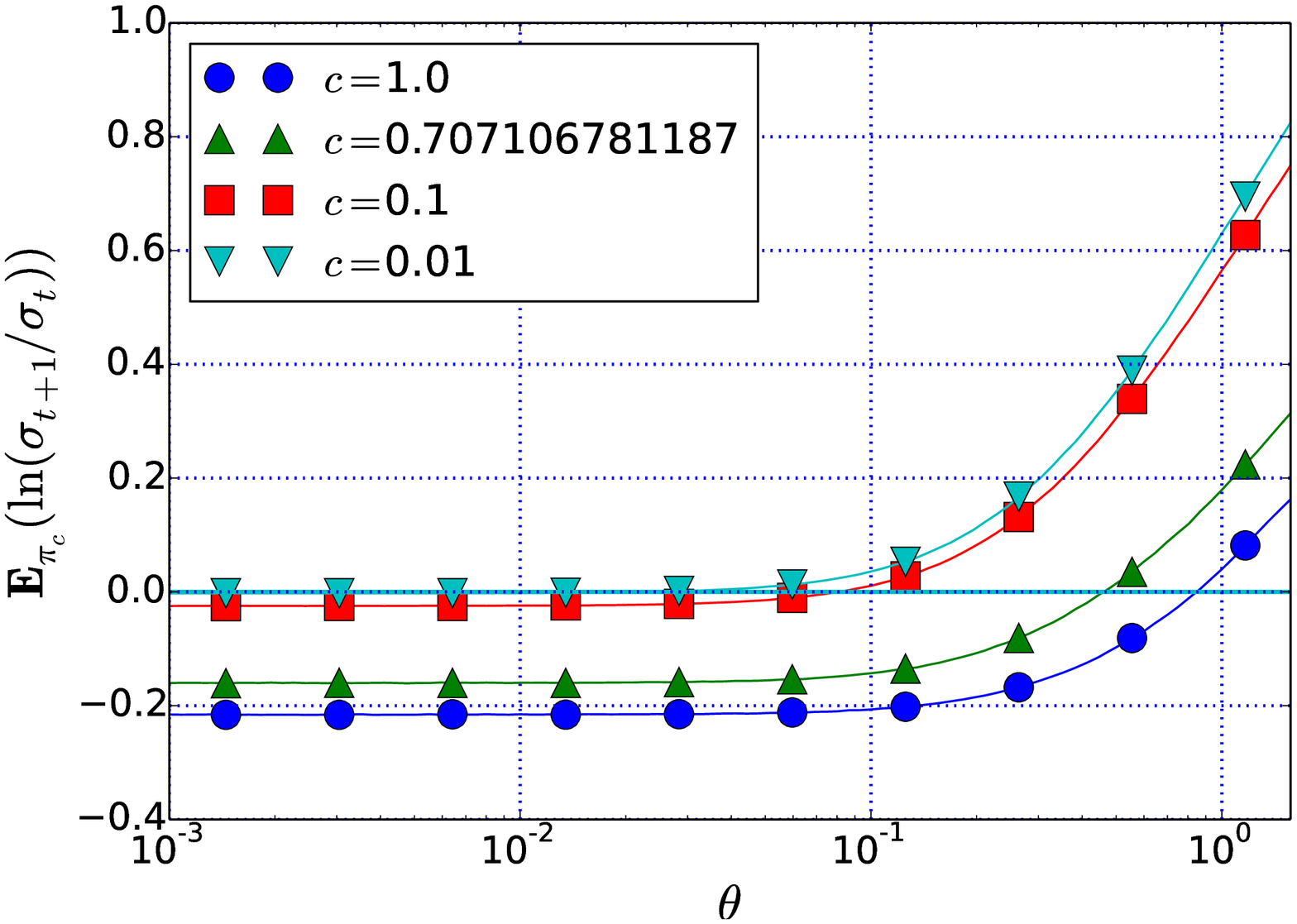}
\end{center}
\caption{Average of the logarithmic adaptation response $\Delta_t = \ln(\sigma_{t+1}/\sigma_t)$ for the $(1,\lambda)$-CSA-ES plotted against the constraint angle $\theta$, for $\lambda = 5$, $c \in \lbrace 1, 1/\sqrt{2}, 0.1, 0.01 \rbrace$\changeaddp{, $d_\sigma = 1$} and dimension $2$. Values below zero (straight line) indicate premature convergence.
}
\label{fg:csasvstcs}
\end{figure}

In Fig.~\ref{fg:csasvstcs} like in the previous Fig.~\ref{fg:csasvst}, the \adp{adaptation response $\Delta_t$ is averaged}\rep{left hand side of Eq.~\eqref{eq:csa} is simulated} for $10^6$ time steps \adp{and plotted} against the constraint angle $\theta$, this time for different values of the cumulation parameter $c$. A lower value of $c$ yields a higher divergence rate for the step-size although $\E_{\pi_c}(\ln(\sigma_{t+1}/\sigma_t))$ appears to converge quickly \new{to an asymptotic constant when $\ln(c)\to-\infty$.}\del{when $c \rightarrow 0$.} Lower values of $c$ hence also allow success of the step-size adaptation for wider range values of $\theta$, and in case of premature convergence a lower value of $c$ means a lower convergence rate.\niko{I would believe that this is somehow an artifact in that $c\ll 1/n$ does not provide adequate convergence properties in general. }

\begin{figure}
\begin{center}
\includegraphics[width=0.49\textwidth]{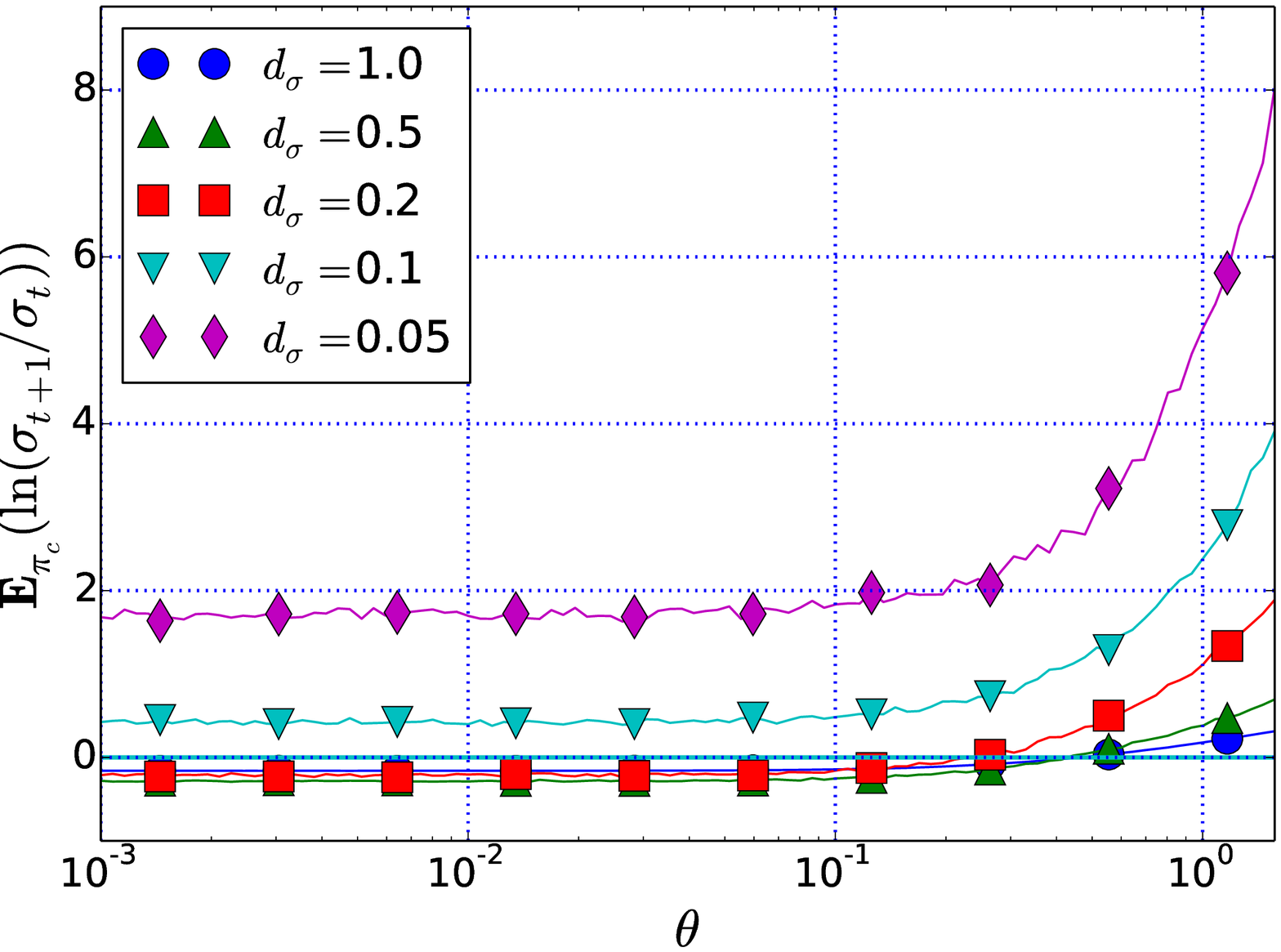}
\end{center}
\caption{Average of the logarithmic adaptation response $\Delta_t = \ln(\sigma_{t+1}/\sigma_t)$ for the $(1,\lambda)$-CSA-ES plotted against the constraint angle $\theta$, for $\lambda = 5$, $c = 1/\sqrt{2}$, $d_\sigma \in \lbrace 1, 0.5, 0.2, 0.1, 0.05 \rbrace$ and dimension 2. Values below zero (straight line) indicate premature convergence.
}
\label{fg:csasvstds}
\end{figure}

\begin{figure}
\begin{center}
\includegraphics[width=0.49\textwidth]{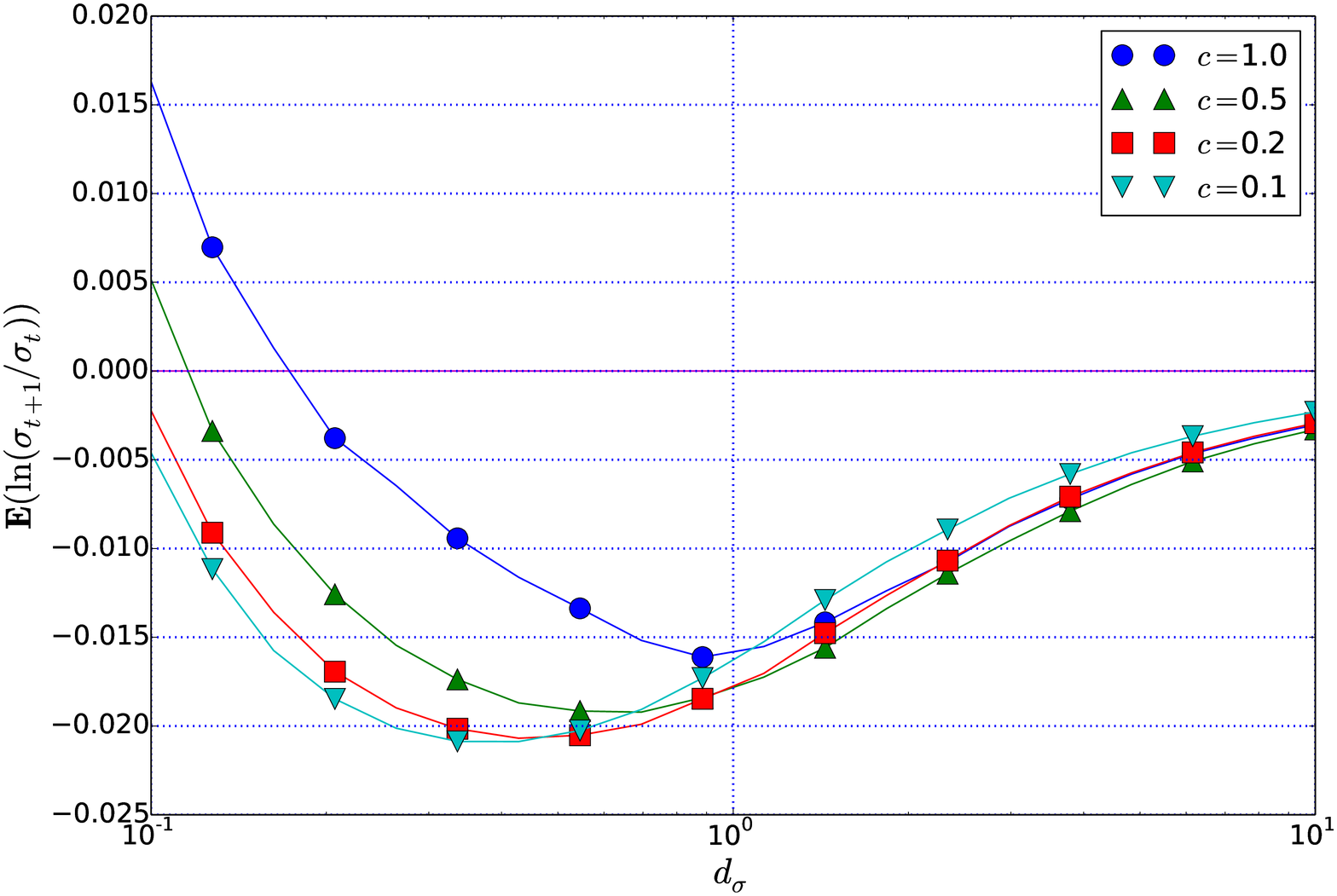}
\end{center}
\caption{Average of the logarithmic adaptation response $\Delta_t = \ln(\sigma_{t+1}/\sigma_t)$ against $d_\sigma$ for the $(1,\lambda)$-CSA-ES minimizing a sphere function for $\lambda = 5$, $c \in \lbrace 1, 0.5, 0.2, 0.1\rbrace$, and dimension \adp{$30$}\rep{$5$}.}
\label{fg:csaspheresvstimeds}
\end{figure}

In Fig.~\ref{fg:csasvstds} the \adp{adaptation response $\Delta_t$ is averaged}\rep{left hand side of Eq.~\eqref{eq:csa} is simulated} for $10^4$ time steps for the $(1,\lambda)$-CSA-ES plotted against the constraint angle $\theta$, for $\lambda = 5$, $c = 1/\sqrt{2}$, $d_\sigma \in \lbrace 1, 0.5, 0.2, 0.1, 0.05 \rbrace$ and dimension $2$.\rep{A lower value of $d_\sigma$  allows larger change of step-size and induces here a bias towards increasing the step-size. This is confirmed in Fig.~\ref{fg:csasvstds} where a} A low enough value of $d_\sigma$ implies geometric divergence \adp{of the step-size} regardless of the constraint angle.  \adp{However, simulations suggest that while for $d_\sigma \geq 1$ the probability that $f(\Xt) > f(\Xt[0])$ is close to $1$, this probability decreases with smaller values of $d_\sigma$.} A low value of $d_\sigma$ will also prevent convergence when it is desired, as shown in Fig.~\ref{fg:csaspheresvstimeds}.

In Fig.~\ref{fg:csaspheresvstimeds} the average of $\ln(\stt/\st)$ is plotted against $d_\sigma$ for the $(1,\lambda)$-CSA-ES  minimizing a sphere function $f_{\textrm{sphere}} : \x \mapsto \| \x \|$, for $\lambda = 5$, $c \in \lbrace 1, 0.5, 0.2, 0.1\rbrace$ and dimension \adp{$30$}\rep{$5$}, averaged over $10$ runs. Low values of $d_\sigma$ \rep{induce a bias towards increasing the step-size, which}make\rep{s} the algorithm diverge while convergence is desired here.


\begin{figure} 
\begin{center}
\includegraphics[width=0.49\textwidth]{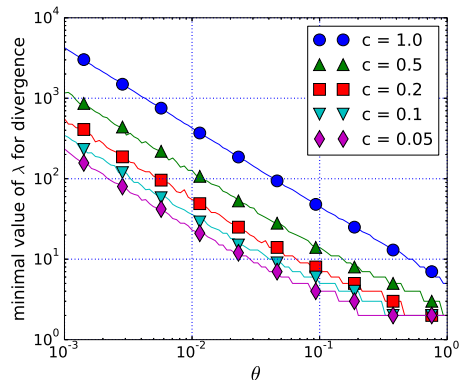}
\end{center}
\caption{Minimal value of $\lambda$ allowing geometric divergence for the $(1,\lambda)$-CSA-ES plotted against the constraint angle $\theta$, for $c \in \lbrace 1., 0.5, 0.2, 0.05 \rbrace$\changeaddp{, $d_\sigma = 1$} and dimension $2$.}
\label{fg:lambdacrit}
\end{figure}

In Fig.~\ref{fg:lambdacrit}, the smallest population size allowing geometric divergence \adp{on the linear constrained function} is plotted against the constraint angle for different values of $c$. Any value of $\lambda$ \changeremp{over}\changeaddp{above} the curve \adp{implies the geometric divergence of the step-size}\rep{will succeed in solving the linear function on this constraint problem} for the corresponding values of $\theta$ and $c$. We see that lower values of $c$ allow for lower values of $\lambda$. It appears that the required value of $\lambda$ scales inversely proportionally with $\theta$. These curves were plotted by simulating runs of the algorithm for different values of $\theta$ and $\lambda$, and stopping the runs when the logarithm of the step-size had decreased or increased by $100$ (for $c=1$) or $20$ (for the other values of $c$). If the step-size had decreased (resp. increased) then this value of $\lambda$ became a lower (resp. upper) bound for $\lambda$ and a \changeremp{higher}\changeaddp{larger} (resp. \changeremp{lower}\changeaddp{smaller}) value of $\lambda$ would be tested until the estimated upper and lower bounds for $\lambda$ would meet. \adp{Also, simulations suggest that for increasing values of $\lambda$ the probability that $f(\Xt) > f(\Xt[0])$ increases to $1$, so large enough values of $\lambda$ appear to solve the linear function on this constrained problem\new{, as expected}.}

\begin{figure} 
\begin{center}
\includegraphics[width=0.49\textwidth]{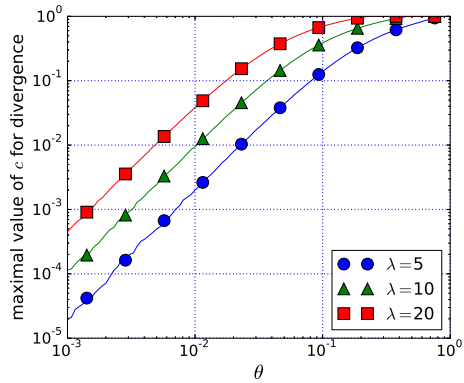}
\end{center}
\caption{Transition boundary for $c$ between convergence and divergence (lower value of $c$ is divergence) for the $(1,\lambda)$-CSA-ES plotted against the constraint angle $\theta$, for $\lambda \in \lbrace 5, 10, 20 \rbrace$ and dimension $2$.
}
\label{fg:ccrit}
\end{figure}

In Fig.~\ref{fg:ccrit} the \new{largest}\del{highest} value of $c$ leading to geometric divergence \adp{of the step-size} is plotted against the constraint angle $\theta$ for different values of $\lambda$. We see that \changeremp{higher}\changeaddp{larger} values of $\lambda$ allow higher values of $c$ to be taken, and when $\theta \rightarrow 0$ the critical value of $c$ appears proportional to $\theta^2$. These curves were plotted following a similar scheme than with Fig.~\ref{fg:lambdacrit}. For a certain $\theta$ the algorithm is ran with a certain value of $c$, and when the logarithm of the step-size has increased (resp. decreased) by more than $1000\sqrt{c}$ the run is stopped, the value of $c$ tested becomes the new lower (resp. upper) bound for $c$ and a new $c$ taken between the lower and upper bounds is tested, until the lower and upper bounds are distant by less than the precision $\theta^2/10$. \adp{Similarly as with $\lambda$, simulations suggest that for small enough values of $c$ the probability that $\lim_{t \to +\infty} f(\Xt) > f(\Xt[0])$ is equal to $1$, so small enough values of $c$ appear to solve the linear function on this constrained problem.} 

\section{Discussion}\label{sc:discuss}

We investigated the $(1,\lambda)$-ES with constant step-size and cumulative step-size adaptation optimizing a linear function under a linear constraint handled by resampling unfeasible solutions. In the case of constant step-size or cumulative step-size adaptation when $c=1$ we prove the stability (formally $V$-geometric ergodicity) of the Markov chain $\dmarkov$ defined as the normalised distance to the constraint, which was \emph{presumed} in \citet{arnold2011behaviour}. This property implies the divergence of the algorithm with constant step-size at a constant speed (see Theorem~\ref{th:cstdiv}), and the geometric divergence or convergence of the algorithm with step-size adaptation (see Theorem~\ref{th:csa}). In addition, it ensures (fast) convergence of Monte Carlo simulations of the divergence rate, justifying their use.

In the case of cumulative step-size adaptation simulations suggest that geometric divergence occurs for a small enough cumulation parameter, $c$, or large enough population size, $\lambda$. 
In simulations we find the critical values \adp{with constraint angle}\rep{for} $\theta\to0$ following $c\propto \theta^2$ \new{or}\del{and} $\lambda\propto 1 / \theta$. Smaller values of the constraint angle seem to increase the difficulty of the problem arbitrarily, i.e.\ no given values for $c$ and $\lambda$ solve the problem for \emph{every} $\theta \in (0,\pi/2)$. 
\adp{However, when using a repair method to handle the constraint instead of resampling with the $(1,\lambda)$-CSA-ES, fixed values of $\lambda$ and $c$ can solve the problem for every $\theta \in (0, \pi/2)$ \cite{arnold2013resamplingvsrepair}.}

Using a different covariance matrix to generate new samples \rep{can be interpreted as}\adp{implies} a change of the constraint angle (see \citealt{ch2014linearproblemgendist} for more details). Therefore, \rep{a correct }adaptation of the covariance matrix \rep{will}\adp{may} render the problem arbitrarily close to the \new{most simple} one with $\theta=\pi/2$. The unconstrained linear function case has been shown to be solved by a $(1,\lambda)$-ES with cumulative step-size adaptation for a population size larger than $3$, regardless of other internal parameters \cite{chotard2012TRcumulative}. 
We believe this is one reason for using covariance matrix adaptation with ES when dealing with constraints\adp{, as has been done in \cite{arnold2012cmaconstraint}}, as pure step-size adaptation has been shown to be liable to fail on even a very basic problem.\niko{In general, I still believe that it is better to use a penalty approach to deal with constraints efficiently. }

This work provides a methodology that can be applied to many ES variants. It demonstrates that a rigorous analysis of the constrained problem can be achieved.  It relies on the theory of Markov chains for a continuous state space that once again proves to be a natural theoretical tool for analyzing ESs, complementing particularly well previous studies \cite{arnold2011behaviour,arnold2012behaviour,arnold2008behaviour}.

\section*{Acknowledgments}This work was supported by the grants  ANR-2010-COSI-002 (SIMINOLE) and ANR-2012-MONU-0009 (NumBBO) of the French National Research Agency.

\small

\bibliographystyle{apalike}
\bibliography{biblio}

\normalsize

\section*{Appendix} \label{sc:appendix}

Proof of Lemma~\ref{lm:farNtstar}.
\begin{myproof}
From Proposition~\ref{pr:g} \adp{and Lemma~\ref{lm:pstar}} the density probability function of $\G(\delta, \W)$ is $\pstar$, and from Eq.~\eqref{eq:pstar}
\begin{equation*}
 \pstar\left(  \left(  \begin{array}{c} x \\ y \end{array}  \right)  \right) = \lambda\frac{\varphi(x)\varphi(y)\1_{\ddomain} \! \left(\delta - \left(  \begin{array}{c} x \\ y \end{array}\right).\n\right)}{\Phi(\delta)} \cdfnone(x)^{\lambda - 1} \enspace \re{.}\ad{,}
\end{equation*}
\adp{where $\cdfnone$ is the cumulative density function of $[\G(\delta, \W)]_1$, whose probability density function is $\ptildeone$.}
From Eq.~\eqref{eq:p1}, $\ptildeone(x) = \varphi(x) \Phi((\delta - x \cos\theta)/\sin\theta)/\Phi(\delta)$, so as $\delta \rep{\geq}\adp{>} 0$ we have $1 \geq \Phi(\delta) \rep{\geq}\adp{>} \Phi(0) = 1/2$, hence $\ptildeone(x) \rep{\leq}\adp{<} 2\varphi(x)$. So $\ptildeone(x)$ converges when $\delta \rightarrow +\infty$ to $\varphi(x)$ while being bounded by $2\varphi(x)$ which is integrable. Therefore we can apply Lebesgue's dominated convergence theorem: $\cdfnone$ converges to $\Phi$ when $\delta \rightarrow + \infty$ and is finite.

For $\delta \in \ddomain$ and $(x,y)\in\R^2$ let $h_{\delta,y}(x)$ be $\exp(a x) \pstar((x,y))$. With Fubini-Tonelli's theorem  $\E(\exp( \G(\delta, \W).(a,b))) = \int_\R \int_\R \exp(b y) h_{\delta,y}(x) \mathrm{d}x \mathrm{d}y$. For $\delta \rightarrow +\infty$, $h_{\delta,y}(x)$ converges to $\exp(a x) \lambda \varphi(x) \varphi(y) \Phi(x)^{\lambda-1}$ while being dominated by $2\lambda \exp (a x)\varphi(x)\varphi(y)$, which is integrable. Therefore by the dominated convergence theorem and as the density of $\Nlambda$ is $x \mapsto \lambda \varphi(x)\Phi(x)^{\lambda-1}$, when $\delta \rightarrow +\infty$, $\int_{\R}h_{\delta,y}(x)\mathrm{d}x$ converges to $ \varphi(y) \E(\exp(a \Nlambda)) < \infty$.

So the function $y \mapsto \exp(b y) \int_{\R}h_{\delta,y}(x)\mathrm{d}x$ converges to $y \mapsto \exp(b y) \varphi(y) \E(\exp (a \Nlambda ))$ while being dominated by $y \mapsto 2\lambda\varphi(y)\exp(by) \int_\R \exp(ax)\varphi(x)\mathrm{d}x$ which is integrable. Therefore we may apply the dominated convergence theorem: $\E( \exp( \G(\delta, \W).(a,b) ) )$ converges to $\int_{\R} \exp(b y) \varphi(y) \E(\exp(a \Nlambda) ) \mathrm{d}y$ which equals to $\E(\exp(a \Nlambda))\E(\exp(b\Nlun))$; and this quantity is finite.

The same reasoning can be applied to $E(\bar{K})$.
\end{myproof}

Proof of Lemma~\ref{lm:farcsa}.
\begin{myproof}
\adp{As in Lemma~\ref{lm:farNtstar}, let $E_1$, $E_2$ and $E_3$ denote respectively $\E( \exp ( -\frac{\alpha}{2d_\sigma n} ( \Nlambda^2 -1 ) ) )$, $\E( \exp ( -\frac{\alpha}{2d_\sigma n} ( \Nl(0,1)^2 -1 ) ) )$, and $\E( \exp ( -\frac{\alpha}{2d_\sigma n} ( K - n + 2 ) ) )$, where $K$ is a random variable following a chi-squared distribution with $n-2$ degrees of freedom.
}
Let us denote $\pdfchi$ the probability density function of $K$. \adp{Since $\pdfchi (z) = (1/2)^{(n-2)/2}/\Gamma((n-2)/2) z^{(n-2)/2} \exp(-z/2)$,  $E_3$ is finite.}

Let $h_\delta$ be a function such that for $(x,y)\in \R^{2}$
\begin{equation*}
h_\delta(x,y) = \frac{\left|\delta -  a x - b y \right|^\alpha} {\rep{\delta^\alpha} \exp\left(\frac{\alpha}{2d_\sigma n}\left(x^2 + y^2 - 2\right)\right)} \enspace ,
\end{equation*}
where $a := \cos \theta$ and $b := \sin \theta$.

\adp{From Proposition~\ref{pr:g} and Lemma~\ref{lm:pstar}, t}\rep{T}he probability density function of $(\G(\delta, \Wt), K)$ is $\pstar \pdfchi$. Using the theorem of Fubini-Tonelli the expected value of the random variable $\frac{\left(\delta -  \G(\delta, \Wt).\n \right)^\alpha} {\rep{\delta^\alpha} \stmod[1](\delta,\W, K)\adp{^\alpha}}$, that we denote $E_\delta$, is
\begin{align*} 
E_\delta  &\adp{= \int_{\R}\int_{\R}\int_{\R} \frac{|\delta - ax -by|^\alpha \pstar ( (x,y) ) \pdfchi(z) }{\exp\left(\frac{\alpha}{2d_\sigma }\left( \frac{\|(x,y)\|^2 + z}{n} - 1 \right)\right)}  \mathrm{d}z \mathrm{d}y \mathrm{d}x}  \\
 &\adp{= \int_{\R}\int_{\R}\int_{\R} \frac{|\delta - ax -by|^\alpha \pstar ( (x,y) ) \pdfchi(z) }{\exp\left(\frac{\alpha}{2d_\sigma n}\left( x^2 + y^2   - 2 \right)\right)\exp\left(\frac{\alpha}{2d_\sigma n }\left(  z - (n-2)\right)\right)}  \mathrm{d}z \mathrm{d}y \mathrm{d}x}  \\
 &=  \int_{\R}\int_{\R}\int_{\R} \frac{h_\delta(x,y) \pstar ( (x,y) ) \pdfchi(z) }{\exp\left(\frac{\alpha}{2d_\sigma \adp{n} }\left(  z - (n-2)\right)\right)}  \mathrm{d}z \mathrm{d}y \mathrm{d}x \enspace .
\end{align*}
Integration over $z$ yields $E_\delta = \int_{\R}\int_{\R} h_\delta(x,y) \pstar ( (x,y) ) \mathrm{d}y \mathrm{d}x E_3$. 

\adp{We now study the limit when $\delta \to +\infty$ of $E_\delta/\delta^\alpha $.}
\newcommand{\pdfnl}{\varphi_{\Nlambda}}
Let $\pdfnl$ denote the probability density function of $\Nlambda$. \adp{For all $\delta \in \R_+^*$, $\Phi(\delta) > 1/2$, and for all $x \in \R$, $\cdfnone(x) \leq 1$, hence with  \eqref{eq:p} and \eqref{eq:pstar}
\begin{equation} \label{eq:pstardom}
\pstar(x,y) = \lambda \frac{\varphi(x)\varphi(y)\1_{\R_+^*}(\delta - ax - by)}{\Phi(\delta)} \cdfnone(x)^{\lambda - 1} \leq \lambda \frac{\varphi(x) \varphi(y)}{\Phi(0)} \enspace ,
\end{equation}
and w}\rep{W}hen $\delta \rightarrow +\infty$, as shown in the proof of Lemma~\ref{lm:farNtstar}, $\pstar((x,y) )$ converges to $\pdfnl(x)\varphi(y)$\adp{.}\rep{while being dominated by $\lambda\varphi(x)\varphi(y)/\Phi(0)$.} For $\delta \geq 1$, $|\delta -ax -by|/\delta \leq 1 + |ax + by|$ with the triangular inequality. Hence 
\adp{
\begin{align} \label{eq:infEd}
 \pstar((x,y)) \frac{h_\delta(x,y)}{\delta^\alpha}  &\leq \lambda\frac{\varphi(x)\varphi(y)}{\Phi(0)} \frac{(1 + |ax + by|)^\alpha }{\exp\left(\frac{\alpha}{2d_\sigma n} \left( x^2 + y^2 - 2 \right) \right) } ~~~~ \textrm{for }\delta \geq 1, \textrm{ and} \\
 \label{eq:infEdcv}
 \pstar((x,y)) \frac{h_\delta(x,y)}{\delta^\alpha} &\underset{\delta \to +\infty}{\longrightarrow} \pdfnl(x)\varphi(y) \frac{1}{\exp\left(\frac{\alpha}{2d_\sigma n} \left( x^2 + y^2 - 2 \right) \right)} \enspace .
\end{align}
Since the right hand side of \eqref{eq:infEd} is integrable, we can use Lebesgue's dominated convergence theorem, and deduce from \eqref{eq:infEdcv} that 
\begin{align*}
\frac{E_\delta}{\delta^\alpha} = \int_{\R}\int_{\R} \frac{h_\delta(x,y)}{\delta^\alpha} \pstar ( (x,y) ) \mathrm{d}y \mathrm{d}x E_3 & \underset{\delta \to +\infty}{\longrightarrow} \int_{\R}\int_{\R}\frac{\pdfnl(x)\varphi(y) }{\exp\left(\frac{\alpha}{2d_\sigma n} \left( x^2 + y^2 - 2 \right) \right)} \mathrm{d}y\mathrm{d}x E_3  \\
\textrm{and so}~~~~ \frac{E_\delta}{\delta^\alpha} & \underset{\delta \to +\infty}{\longrightarrow} E_1 E_2 E_3 < \infty \enspace .
\end{align*}
Since $\delta^\alpha / (\delta^\alpha + \delta^{-\alpha})$ converges to $1$ when $\delta \to +\infty$, $E_\delta/(\delta^\alpha + \delta^{-\alpha})$ converges to $E_1 E_2 E_3$ when $\delta \to + \infty$.
}

\rep{
$h_\delta(x,y) \pstar((x,y) ) $ converges to $\exp(-\frac{\alpha}{2d_\sigma n}(x^2 + y^2 - 2)) \pdfnl(x) \varphi(y)$, while being dominated by $(1 + |ax + by|)^\alpha \exp(\frac{- \alpha}{2d_\sigma n}(x^2 + y^2 - 2)) \lambda\varphi(x)\varphi(y)/\Phi(0)$, which is integrable over $y$.
 Therefore, applying Lebesgue's Theorem of dominated convergence yields that $\int_\R h_\delta(x,y) \pstar((x,y) ) \mathrm{d}y$ converges to $\int_\R \exp(-\frac{\alpha}{2d_\sigma n}(x^2 + y^2 - 2)) \pdfnl(x) \varphi(y) \mathrm{d}y$, which equals $ \exp(-\frac{\alpha}{2d_\sigma n}(x^2 - 1)) p_{\Nlambda}(x)E_2$ and this quantity is finite.

Finally, as $\int_\R h_\delta(x,y) \pstar((x,y)) \mathrm{d}y$ is dominated by $\int_\R (1 + |ax + by|)^\alpha \exp(\frac{- \alpha}{2d_\sigma n}(x^2 + y^2 - 2)) \lambda\varphi(x)\varphi(y)/\Phi(0) \mathrm{d}y$ which is integrable, with Lebesgue's dominated convergence theorem we get that $\int_\R\int_\R h_\delta(x,y) \pstar((x,y) ) \mathrm{d}y\mathrm{d}y$ converges to $\int_{\R}\exp(-\frac{\alpha}{2d_\sigma n}(x^2 - 1)) \pdfnl(x)E_2 \mathrm{d}x$, which equals $E_1 E_2$. Therefore $E_\delta$ converges to  $E_1 E_2 E_3$.}

\adp{
We now study the limit when $\delta \to 0$ of $\delta^\alpha E_\delta$, and restrict $\delta$ to $(0,1]$. When $\delta \to 0$, $\delta^\alpha h_\delta(x,y) \pstar((x,y))$ converges to $0$. Since we took $\delta \leq 1$, $|\delta + ax + by| \leq 1 + |ax + by|$, and with \eqref{eq:pstardom} we have
\begin{equation} \label{eq:zEd}
\delta^\alpha h_\delta(x,y) \pstar((x,y)) \leq \lambda \frac{(1+|ax + by|)^\alpha \varphi(x) \varphi(y)}{\Phi(0) \exp\left(\frac{\alpha}{2d_\sigma n}\left( x^2 + y^2 -2 \right) \right)} ~~~~ \textrm{for }0<\delta\leq 1 \enspace .
\end{equation}
The right hand side of \eqref{eq:zEd} is integrable, so we can apply Lebesgue's dominated convergence theorem, which shows that $\delta^\alpha E_\delta$ converges to $0$ when $\delta \to 0$. And since $(1/\delta^\alpha)/(\delta^\alpha+\delta^{-\alpha})$ converges to $1$ when $\delta \to 0$, $E_\delta / (\delta^\alpha + \delta^{-\alpha})$ also converges to $0$ when $\delta \to 0$.
}

\adp{
Let $H_3$ denote $\E(\exp(\alpha/(2d_\sigma n)(K - (n+2))))$. Since $\pdfchi (z) = (1/2)^{(n-2)/2}/\Gamma((n-2)/2) z^{(n-2)/2} \exp(-z/2)$, when $\alpha$ is close enough to $0$, $H_3$ is finite.
Let $H_\delta$ denote the expected value of the random variable $\frac { \stmod[1](\delta,\W, K)^\alpha}{\left(\delta -  \G(\delta, \Wt).\n \right)^\alpha}$, then
\begin{equation*}
H_\delta = \int_{\R}\int_{\R}\int_{\R} \frac{  \pstar ( (x,y) ) \pdfchi(z) \exp\left(\frac{\alpha}{2d_\sigma n}\left( z - (n-2)\right)\right) }{h_\delta(x,y) }  \mathrm{d}z \mathrm{d}y \mathrm{d}x \enspace .
\end{equation*}
Integrating over $z$ yields $H_\delta = \int_{\R}\int_{\R} \frac{\pstar ( (x,y) )  }{h_\delta(x,y)  }   \mathrm{d}y \mathrm{d}x H_3$.
}

\adp{
We now study the limit when $\delta \to +\infty$ 
 of $H_\delta / \delta^\alpha$. 
 With \eqref{eq:pstardom}, we have that
\begin{align*} \label{eq:Hd}
\frac{\pstar((x,y))}{ \delta^\alpha h_\delta(x,y)} &\leq  \lambda\frac{\varphi(x)\varphi(y)}{\Phi(0)}\frac{ \exp \left(\frac{\alpha}{2d_\sigma n}\left( x^2 + y^2 - 2\right) \right)}{ \delta^\alpha |\delta - a x  - b y|^\alpha } \enspace .
\end{align*}
With the change of variables $\tilde{x} = x - \delta/a$ we get
\begin{align*}
\frac{\pstar((\tilde{x}+\frac{\delta}{a},y))}{ \delta^\alpha h_\delta(\tilde{x}+\frac{\delta}{a},y)} &\leq  \lambda\frac{\exp\left(-\frac{(\tilde{x}+\frac{\delta}{a})^2}{2}\right) \varphi(y)}{\sqrt{2\pi}\Phi(0)} \frac{ \exp \left(\frac{\alpha}{2d_\sigma n}\left(\left(\tilde{x}+\frac{\delta}{a}\right)^2 + y^2 - 2\right) \right)}{ \delta^\alpha |a \tilde{x}  + b y|^\alpha }  \\
 &\leq  \lambda\frac{\varphi(\tilde{x}) \varphi(y)}{\Phi(0)}  \frac{\exp\left( \frac{\alpha}{2d_\sigma n}\left( \tilde{x}^2 + y^2 - 2 \right)\right)}{|a\tilde{x} + by|^\alpha}    \frac{ \exp\left( \! \left( \! \frac{\alpha}{2d_\sigma n} - \frac{1}{2}\right)\left(2\frac{\delta}{a}\tilde{x} + \frac{\delta^2}{a^2} \! \right) \! \right)}{\delta^\alpha} \\
 &\leq \lambda\frac{\varphi(\tilde{x}) \varphi(y)}{\Phi(0)} \frac{1}{h_0(\tilde{x}, y)} \frac{ \exp\left( \left( \frac{\alpha}{2d_\sigma n} - \frac{1}{2}\right)\left(2\frac{\delta}{a}\tilde{x} + \frac{\delta^2}{a^2}\right) \right) } {\exp(\alpha \ln(\delta))} \enspace .
\end{align*}
}
\adp{
An upper bound for all $\delta \in \R_+^*$ of the right hand side of the previous inequation \del{corresponds to an upper bound to the function $l : \delta \in \R_+^* \mapsto (\alpha/(2d_\sigma n) - 1/2)(2(\delta/a) \tilde{x} + \delta^2/a^2) - \alpha \ln(\delta)$\todo{to the exponential of the upper bound of this function}.} is a function of an upper bound of the function $l : \delta \in \R_+^* \mapsto (\alpha/(2d_\sigma n) - 1/2)(2(\delta/a) \tilde{x} + \delta^2/a^2) - \alpha \ln(\delta)$.  And since we are interested in a limit when $\delta \to + \infty$, we can restrict our search of an upper bound of $l$ to $\delta \geq 1$. Let $c:= \alpha/(2d_\sigma n) - 1/2$. We take $\alpha$ small enough to ensure that $c$ is negative. An upper bound to $l$ can be found through derivation:
\begin{align*}
\frac{\partial l(\delta)}{\partial \delta} = 0 &\Leftrightarrow 2\frac{c}{a^2}\delta + 2 \frac{c}{a}\tilde{x} - \frac{\alpha}{\delta} = 0 \\
&\Leftrightarrow 2\frac{c}{a^2}\delta^2 + 2 \frac{c}{a}\tilde{x}\delta - \alpha = 0
\end{align*}
The discriminant of the quadratic equation is $\Delta = 4(c^2/a^2) \tilde{x}^2 + 8\alpha c/a^2$. \del{Since we restricted $\delta$ to $[1, +\infty)$ and that the derivative of $l$ is a quadratic function with a negative quadratic coefficient $2c/a^2$, \todo{delta times the derivative of l is a quadratic function (it changes a bit the reasoning, you have to argue with the maximum of a continuous function I guess)}} The derivative of $l$ multiplied by $\delta$ is a quadratic function   with a negative quadratic coefficient $2c/a^2$. Since we restricted $\delta$ to $[1, +\infty)$, multiplying the derivative of $l$ by $\delta$ leaves its sign unchanged. So the maximum of $l$ is attained for $\delta$ equal to $1$ or for $\delta$ equal to $\delta_M := (-2c/a\tilde{x} - \sqrt{\Delta})/(4c/a^2)$, and so $l(\delta) \leq \max(l(1), l(\delta_M))$ for all $\delta \in [1, +\infty)$. We also have that $\lim_{\tilde{x} \to \infty }\sqrt{\Delta}/\tilde{x} = 2|c|/a = -2ca$, so $\lim_{\tilde{x} \to \infty} \delta_M / \tilde{x} = (-2c/a - (-2c/a))/(4c/a^2) = 0$. Hence when $|\tilde{x}|$ is large enough, $\delta_M \leq 1$, so since we restricted $\delta$ to $[1,+\infty)$ there exists $m > 0$ such that if $|\tilde{x}| > m$, $l(\delta) \leq l(1)$ for all $\delta \in [1, +\infty)$. And trivially, $l(\delta)$ is bounded for all $\tilde{x}$ in the compact set $[-m,m]$ by a constant $M> 0$, so $l(\delta) \leq \max(M, l(1)) \leq M + |l(1)|$ for all $\tilde{x} \in \R$ and all $\delta \in [1, +\infty)$. Therefore
\begin{align*}
\frac{\pstar((\tilde{x}+\frac{\delta}{a},y))}{ \delta^\alpha h_\delta(\tilde{x}+\frac{\delta}{a},y)} &\leq \lambda\frac{\varphi(\tilde{x}) \varphi(y)}{\Phi(0)} \frac{1}{h_0(\tilde{x}, y)} \exp(M + |l(1)|) \\
 &\leq 
\lambda\frac{\varphi(\tilde{x}) \varphi(y)}{\Phi(0)} \frac{1}{h_0(\tilde{x}, y)} \exp\left(M + \left|2\frac{c}{a}\tilde{x} + \frac{c}{a^2}\right| \right)
 \enspace .
\end{align*}
For $\alpha$ small enough, the right hand side of the previous inequation is integrable. And since the left hand side of this inequation converges to $0$ when $\delta \to +\infty$, according to Lebesgue's dominated convergence theorem  $H_\delta / \delta^\alpha $ converges to $0$  when $\delta \to +\infty$. And since $\delta^\alpha / (\delta^\alpha + \delta^{-\alpha})$ converges to $1$ when $\delta \to +\infty$, $H_\delta / (\delta^\alpha + \delta^{-\alpha})$ also converges to $0$ when $\delta \to +\infty$.
}

\adp{
We now study the limit when $\delta \to 0$ of $H_\delta/(\delta^\alpha + \delta^{-\alpha})$. Since we are interested in the limit for $\delta \to 0$, we restrict $\delta$ to $(0,1]$. Similarly as what was done previously, with the change of variables $\tilde{x} = x - \delta/a$, 
\begin{align*}
\frac{\pstar((\tilde{x}+\frac{\delta}{a},y))}{(\delta^\alpha + \delta^{-\alpha})h_\delta(\tilde{x}+\frac{\delta}{a},y)} &\leq \lambda\frac{\varphi(\tilde{x}) \varphi(y)}{\Phi(0)} \frac{1}{h_0(\tilde{x}, y)} \frac{\exp\left( \left( \frac{\alpha}{2d_\sigma n} - \frac{1}{2}\right)\left(2\frac{\delta}{a}\tilde{x} + \frac{\delta^2}{a^2}\right) \right)}{\delta^\alpha + \delta^{-\alpha}} \\
 &\leq \lambda\frac{\varphi(\tilde{x}) \varphi(y)}{\Phi(0)h_0(\tilde{x}, y)} \exp\left( \left( \frac{\alpha}{2d_\sigma n} - \frac{1}{2}\right)\left(2\frac{\delta}{a}\tilde{x} + \frac{\delta^2}{a^2}\right) \right)   \enspace . 
\end{align*}
Take $\alpha$ small enough to ensure that $\alpha/(2d_\sigma n) - 1/2$ is negative. \del{Then an upper bound for $\delta \in (0,1]$ of the right hand side of the previous inequality corresponds to the maximum of the function $k : \delta \in (0,1] \mapsto 2\delta \tilde{x}/a + \delta^2 /a^2$}\todo{here again it's litterally not what you want to say, the upper bound of $k$ will not be an upped bound of the previous inequality.} Then an upper bound for $\delta \in (0,1]$ of the right hand side of the previous inequality is a function of an upper bound of the function $k : \delta \in (0,1] \mapsto 2\delta \tilde{x}/a + \delta^2 /a^2$. This upper bound can be found through derivation: $\partial k (\delta)/\partial \delta = 0$ is equivalent to $2\tilde{x}/a + 2\delta /a^2 = 0$, and so the upper bound of $k$ is realised at $\delta_M := -a \tilde{x}$. However, since we restricted $\delta$ to $(0,1]$, for $\tilde{x} \geq 0$ we have $\delta_M \leq 0$ so an upper bound of $k$ in $(0,1]$ is realized at $0$, and for $\tilde{x} \leq -1/a$ we have $\delta_M \geq 1$ so  the maximum of $k$ in $(0,1]$ is realized at $1$. Furthermore, $k(\delta_M) = -2\tilde{x}^2 + \tilde{x}^2 = -\tilde{x}^2$ so when $-1/a < \tilde{x} < 0$, $  k(\delta) < 1/a^2$. Therefore $k(\delta) \leq \max(k(0), k(1), 1/a^2)$. Note that $k(0) = 0$ which is inferior to $1/a^2$, and note that $k(1) = 2c\tilde{x}/a +/a^2$. Hence $k(\delta) \leq \max( 2\tilde{x}/a +1/a^2, 1/a^2) \leq |2 \tilde{x} /a + 1/a^2| + 1/a^2 $, and so 
$$
\frac{\pstar((\tilde{x}+\frac{\delta}{a},y))}{(\delta^\alpha + \delta^{-\alpha})h_\delta(\tilde{x}+\frac{\delta}{a},y)} \leq \lambda\frac{\varphi(\tilde{x}) \varphi(y)}{\Phi(0)h_0(\tilde{x}, y)} \exp\left( \left(\frac{\alpha}{2d_\sigma n} - \frac{1}{2}\right) \left( \left|2\frac{\tilde{x}}{a} + \frac{1}{a^2}\right| + \frac{1}{a^2} \right) \right)  \enspace . 
$$
For $\alpha$ small enough the right hand side of the previous inequation is integrable. Since the left hand side of this inequation converges to $0$ when $\delta \to 0$, we can apply Lebesgue's dominated convergence theorem, which proves that $H_\delta/(\delta^\alpha + \delta^{-\alpha})$ converges to $0$ when $\delta \to 0$.
}

\end{myproof}

\end{document}